\documentclass[12pt]{amsart} 
\usepackage{amssymb}
\usepackage{latexsym}
\usepackage{bbm}
\usepackage{amsmath}
\usepackage{amsthm}
\usepackage{amsfonts}

\usepackage{mathabx}
\usepackage{enumerate}
\usepackage{mathrsfs}

\newtheorem{theorem}{Theorem}[section]
\newtheorem{lemma}[theorem]{Lemma}
\newtheorem{corollary}[theorem]{Corollary}
\newtheorem{proposition}[theorem]{Proposition}
\theoremstyle{definition}
\newtheorem{remark}[theorem]{Remark}

\numberwithin{equation}{section}

\newcommand{\C}{{\mathbb C}}
\newcommand{\R}{{\mathbb R}}
\newcommand{\N}{{\mathbb N}}

\newcommand{\sa}{\sigma}

\newcommand{\Bb}{\mathscr{B}}
\newcommand{\Rr}{\mathscr{R}}
\newcommand{\Normx}{\Vert \cdot \Vert_X}
\newcommand{\fn}{\{f_n\}_{n=1}^\infty}
\newcommand{\eps}{\varepsilon}

\newcommand{\lbal}{\underline{\alpha}}
\newcommand{\ubal}{\overline{\alpha}}
\newcommand{\si}{\mbox{sim}}

\newcommand{\sub}{\subseteq}
\newcommand{\whX}{\widehat{X}}

\begin{document}

\title[The finite Hilbert transform]
{Extension and Integral Representation of the finite Hilbert 
Transform In Rearrangement Invariant Spaces}

\dedicatory{Dedicated to the memory of Joe Diestel}




\author[G. P. Curbera]{G. P. Curbera}\address{Guillermo P. Curbera, 
Facultad de Matem\'aticas \& IMUS, Universidad de Sevilla, Calle Tarfia s/n, Sevilla 41012, Spain, curbera@us.es} 

\author[S. Okada]{S. Okada}\address{Susumu Okada, School of 
Mathemathics and Physics, University of Tasmania, 
Private Bag 37, Hobart, TAS 7001, Australia, susbobby@grapevine.com.au} 

\author[W. J. Ricker]{W. J. Ricker}\address{Werner J. Ricker, Math.-Geogr.\  Fakult\"at, 
Katholische Universit\"at Eichst\"att-Ingolstadt, 85072 Eichst\"att, Germany, werner.ricker@ku.de}

\thanks{The first author acknowledges the support  of 
MTM2015-65888-C4-1-P, MINECO (Spain).}

\begin{abstract}
The finite Hilbert transform $T$ is a classical (singular) kernel operator which 
is continuous in every rearrangement invariant space $X$ over $(-1,1)$ having non-trivial Boyd indices.
For $X=L^p$, $1<p<\infty$, this operator has been intensively investigated 
since the 1940's (also under the guise of the ``airfoil equation''). Recently, 
the extension and inversion of $T\colon X\to X$ for more general $X$ 
has been studied in \cite{COR}, where it is shown that there exists a larger 
space $[T,X]$, optimal in a well defined sense, which contains $X$ 
continuously and such that $T$ can be extended to a continuous linear 
operator $T\colon[T,X]\to X$. The purpose of this paper is to continue 
this investigation of $T$ via a consideration of the $X$-valued vector measure 
$m_X\colon A\mapsto T(\chi_A)$ induced by $T$ and its associated integration operator
$f\mapsto \int_{-1}^1f\,dm_X$. In particular, we present integral representations of $T\colon X\to X$
based on the $L^1$-space of $m_X$ and other related spaces of integrable functions.
\end{abstract}
\maketitle


\section{Introduction}\label{section1}

The finite Hilbert transform $T$ on $(-1,1)$ is a classical singular integral operator 
with kernel  $(x,t)\mapsto (x-t)^{-1}$. This is a rather singular kernel with no positivity or monotonicity properties. 
Nevertheless, $T$ maps $L^1(-1,1)$ into weak-$L^1(-1,1)$ and may be interpreted 
as a ``truncated '' Hilbert transform.  For the definition of the Hilbert transform 
$H$ on $\R$ see \cite[Definition III.4.1]{BS}. With the identification 
$L^1(-1,1)=\left\{\phi\chi_{(-1,1)}:\phi\in L^1_{\mbox{\tiny{loc}}}(\R)\right\}$ we have
\begin{align*}
T(f)=-\chi_{(-1,1)}(H|_{L^1(-1,1)})(f),\hspace{5mm} f\in L^1(-1,1);
\end{align*}
see Section \ref{section2}. The singular integral equation $T(f)=g$ (with $g$ given) 
is known as the airfoil equation, \cite[\S  4-3]{T}, according to its application to 
aerodynamics. Other areas for applications include signal processing, tomography and elasticity problems.

In this paper we deal with $T$ acting on a rearrangement invariant space $X$ over $(-1,1)$ 
with non-trivial lower and upper Boyd indices (i.e., $0<\lbal_X\leq \ubal_X <1$); necessarily 
$X\sub L^1(-1,1)$. For such spaces $X$ the operator $T$ maps $X$ into $X$ and the resulting 
linear operator $T_X:X\to X$ is bounded, \cite[\S 3]{COR}. This generalizes the known result for 
the case $X:= L^p(-1,1)$ with $1<p<\infty$. In order to investigate kernel operators $S$ on a 
Banach function space $Y$ (briefly, B.f.s.), a useful approach is to seek a larger B.f.s.\ $Z$ 
(with certain desirable properties) into which $Y$ is continuously embedded and to which 
$S$ admits a $Y$-valued, continuous  linear extension; see, for example, 
\cite{CR-N},\cite{CR},\cite{CR-Survey}, \cite{ORS}, and the references therein. 
This approach, for $S$ being the operator $T_X:X\to X$, was adopted in \cite{COR}. 
For instance, Theorem 4.6 in \cite{COR} constructs such a larger B.f.s.\ 
$Z=[T,X]$ within $L^1(-1,1)$; see (\ref{2.4}) below for the description of $[T,X]$. 
We improve on this slightly in Proposition \ref{proposition2.7}, where it is shown 
that $[T,X]$ is actually the largest B.f.s.\ within $L^0(-1,1)$ having the same properties. 
Accordingly, $[T,X]$ coincides with the definition given at the beginning of \S 3 in \cite{CR-N}. 
The rest of Section \ref{section2} provides the basic tools, mostly from \cite{COR}, on which 
the arguments in Sections \ref{section3} and \ref{section4} are built.

The construction of $[T,X]$ in \cite{COR} is of an analytic nature. In Section \ref{section3} 
we show that $[T,X]$ can also be interpreted as a space of (scalarly) integrable functions 
with respect to an $X$-valued \emph{vector measure} $m_X$, namely that given by 
$A\mapsto T_X(\chi_A)$ on the Borel subsets of $(-1,1)$. This interpretation allows 
us to investigate $T_X$ and $[T,X]$ from a different perspective.  To do so first 
requires a knowledge of the properties of $m_X$. These are presented in 
Proposition \ref{proposition3.2} and include the fact that the range of $m_X$ is 
neither relatively compact nor  order bounded. Consequently, the operator $T_X$ 
is neither completely continuous nor order bounded (cf. Corollary \ref{corollary3.4}).

Let $L^1(m_X)$ be the B.f.s.\ of all $\C$-valued, $m_X$-integrable functions 
and $L^1_w(m_X)$ that of all $\C$-valued, scalarly $m_X$-integrable functions, 
in which case $L^1(m_X)$ is a closed subspace of $L^1_w(m_X)$; see 
Section \ref{section3}. With $X_a$ denoting the order continuous part of $X$, 
the B.f.s.\ $L^1(m_X)$ is the \emph{largest} B.f.s.\ with order continuous norm, 
into which $X_a$ is continuously embedded and to which the restriction 
$T|_{X_a}:X_a\to X$ admits a continuous, linear, $X$-valued extension. 
The B.f.s.\  $L^1(m_X)$ may not have the Fatou property. On the other 
hand, $L^1_w(m_X)$ is the minimal B.f.s.\ with the Fatou property and 
containing $L^1(m_X)$; see Lemma \ref{lemma3.9}. Our main result 
(cf. Theorem \ref{theorem3.10}) establishes the relationships 
\begin{equation}\label{1.1}
X_a\sub L^1(m_X)\sub [T,X]=L^1_w(m_X),
\end{equation}
in which the first inclusion is a continuous embedding, the second an isometric 
embedding and the identity $[T,X]=L^1_w(m_X)$ with equal norms. Moreover, 
the integration operator $I_{m_X}:f\mapsto \int_{-1}^1 f\; d m_X$ on $L^1(m_X)$ 
coincides with $T|_{L^1(m_X)}$ and provides an integral representation of 
$T|_{X_a}$. Equivalent conditions for a function $f\in [T,X]$ to be actually 
$m_X$-integrable are also provided in Theorem \ref{theorem3.10}. 
In Corollary \ref{corollary3.11} we consider extra conditions under which 
sharper results than those in (\ref{1.1}) are obtained.

A crucial fact  in (\ref{1.1}) is the equality $[T,X]=L^1_w(m_X)$. Whereas the 
integration operator $I_{m_X}$ is an $X$-valued extension of $T|_{X_a}$, it 
typically has no further extension to $L^1_w(m_X) \supseteq L^1(m_X)$; this 
is explained in Section \ref{section4}. However, it turns out that there do exist 
some classical function spaces $\widehat{X}$ over $(-1,1)$ which contain $X$ 
continuously and such that $m_X$, when interpreted as an $\widehat{X}$-valued 
vector measure (denoted by $m_{\widehat{X}}$), has the property that 
$[T,X]\sub L^1(m_{\widehat{X}})$ continuously and the operator $T:[T,X]\to X$ has the integral representation 
\begin{equation}\label{1.2}
T(f)= \int_{-1}^1 f \; d m_{\widehat{X}},\hspace{5mm} f\in [T,X].
\end{equation}
That is, $m_{\widehat{X}}$ is $\widehat{X}$-valued, but the integral 
$\int_{-1}^1f\; d m_{\widehat{X}}$ (typically an element of $\widehat{X}$ 
if $f\in L^1( m_{\widehat{X}})$) actually belongs to $X\sub \widehat{X}$ 
whenever $f$  belongs to the subspace $[T,X]\sub L^1(m_{\widehat{X}})$. 
We point out that the containment  $[T,X]\subseteq L^1(m_{\widehat{X}})$ may be \textit{proper};
see Remark \ref{remark4.2}(ii) for $X=L^{2,r}$ with $2<r\le\infty$, in which case $p_X=2\notin J_X=(1,2)$
and  $\widehat{X}=\cap_{1<p<2}\ L^p$.
Relevant spaces $\widehat{X}$ and the representation (\ref{1.2}) 
are the topic of Section \ref{section4}; see Propositions \ref{proposition4.3} and \ref{proposition4.4}.


\section{Preliminaries and Basic results}\label{section2}

We consider the Lebesgue measure $\mu$ defined on the Borel $\sa$-algebra $\Bb$ of the 
open interval $(-1,1)$. The space of all $\C$-valued, $\Bb$-simple functions on $(-1,1)$ is 
denoted by $\mbox{sim}\Bb$. Let $L^0=L^0(-1,1) ~(=L^0(\mu))$ denote the space of all 
$\C$-valued, $\Bb$-measurable  functions. We identify those functions in $L^0$ which 
coincide a.e. Moreover, $L^0$ is a (complex) vector lattice with respect to the a.e.\  pointwise order. 
An order ideal $X$ of $L^0$ is called a Banach function space, briefly B.f.s., based on the 
measure space $((-1,1),\Bb,\mu)$ if $\mbox{sim}\Bb\subseteq X$ and if $X$ is equipped with 
a lattice norm $\Vert \cdot \Vert_X$ for which $X$ is complete. In this case we speak of 
a B.f.s.\ $X$ over $(-1,1)$. A typical B.f.s.\ over $(-1,1)$ is the usual Lebesgue space 
$L^p(-1,1)=L^p(\mu)$ for $1\leq p\leq \infty$, which we shall denote simply by $L^p$.

A B.f.s.\ $X$ over $(-1,1)$ is said to have the Fatou property if, whenever $\{f_n\}_{n=1}^\infty$ 
is a norm-bounded, increasing sequence of positive functions in $X$ such that $f= \sup_{n\in\N}f_n$ 
in the order of $L^0$, it necessarily follows that $f\in X$ and $\lim_{n\to \infty}\Vert f_n\Vert_X=\Vert f\Vert_X$. 
Since this definition is via increasing sequences $\{f_n\}_{n=1}^\infty$ in $X$, a more natural 
terminology would be the $\sa$-Fatou property or the sequential Fatou property. However, in 
our setting, the definition via increasing \emph{sequences} is equivalent to that by 
increasing \emph{nets}, so that we do not have to distinguish between these two notions. 
Indeed, our B.f.s.\ $X$ is always super Dedekind complete  and hence, order separable  
because $L^0$ is order separable (see \cite[Example 23.3(iv)]{LZ}). Moreover, $X$ has a 
countable order basis because $\mu$ being a finite measure ensures that $X$ has a 
countable maximal disjoint system, which is then necessarily a countable order basis 
(see \cite[Definition 28.4 and Theorem 28.5]{LZ}). So Theorem 113.2 in \cite{ZII} 
assures us that the two definitions are equivalent for a B.f.s.\  over $(-1,1)$.

Note that B.f.s.' in \cite{BS}, by definition, have the Fatou property; see (P3) in 
Definitions I.1.1 and I.1.3 and Lemma I.1.5(i) in \cite{BS}. So, whenever we refer 
to \cite{BS} for B.f.s.' $X$ which are \emph{not} assumed to have the Fatou property, 
we shall clearly state that the corresponding  results in \cite{BS} do \emph{not} require the Fatou property.

Let $X$ be a B.f.s.\ over $(-1,1)$. The norm $\Normx$ is $\sa$-order continuous (resp.\  
order continuous, briefly o.c.) if $\Vert f_n\Vert_X \downarrow 0$ for every sequence 
$\fn\subseteq X$ satisfying $f_n\downarrow 0$ (resp. $\Vert f_\tau\Vert_X\downarrow 0 $ 
for every net $\{f_\tau\}_\tau \subseteq X$ satisfying $f_\tau \downarrow 0$). 
Since $X$ is super Dedekind complete, the norm $\Normx$ is $\sa$-o.c.\  if 
and only if it is o.c. (see, for example, \cite[p.25]{ORS}). 
So, we shall speak of o.c. norms from now on. The terminology absolutely 
continuous norm (briefly, a.c.-norm) is also used, \cite{BS}.

Define the order continuous part of $X$ by 
\begin{align*}
X_a:= \left\{f\in X:|f|\geq f_n\downarrow 0\mbox{ with } f_n\in X \mbox{ implies } \Vert f_n\Vert_X\downarrow 0\right\},
\end{align*}
which is also called the absolutely continuous part of $X$. It follows from \cite[Theorem 102.8]{ZII} 
that $X_a$ is a closed order ideal of $X$. In particular, the norm of $X$ restricted to $X_a$ is o.c.
An important fact is that a function $f\in X$ belongs to $X_a$ if and only if 
$\Vert f \chi_{A_n}\Vert_X\downarrow 0$ for every sequence $\{A_n\}_{n=1}^\infty \subseteq \Bb$ 
satisfying $A_n\downarrow \emptyset $; see Propositions I.3.2 and I.3.5 in \cite{BS}; 
the proofs do \emph{not} require the  Fatou property. 

Let $X_b$ be the closure $\overline{\mbox{sim}\Bb}$ of $\mbox{sim}\Bb$ in $X$. Then $X_b$ 
is a closed order ideal of $X$ such that $X_a\subseteq X_b$; see Theorem I.3.11 of \cite{BS}, 
where a careful examination of its proof reveals that it does \emph{not} require the 
Fatou property. Moreover, $X_a=X_b$ if and only if $\mbox{sim}\Bb \subseteq X_a$. 
For this fact, see Theorem I.3.13 of \cite{BS}; the proof does \emph{not} require  the Fatou property.

Let $X$ be a B.f.s.\ over $(-1,1)$. Its (topological) dual space is denoted by $X^*$ and 
we write $\langle f, x^*\rangle:= x^*(f)$ for $f\in X$ and $x^* \in X^*$. The associate 
space $X'$ of $X$ is defined as the vector space of all $g\in L^0$ such that $f g\in L^1$ 
for every $f\in X$; it is also a B.f.s.\ over $(-1,1)$ with respect to the lattice norm
\begin{align*}
\Vert g\Vert_{X'}:= \sup \left\{\left|\int_{-1}^1f g\; d \mu\right|:f\in X, \Vert f\Vert_X\leq 1\right\}.
\end{align*}
Then $X'$ is a closed linear subspace of $X^*$ as each $g\in X'$ defines the continuous 
linear functional  $f\mapsto \int_{-1}^1 f g\; d \mu$ on $X$, which enables us to write 
$\langle f,g\rangle =\int_{-1}^1 f g \;d\mu$ for $f\in X$. It is also the case that 
$\Vert g \Vert_{X'}=\sup \left\{\int_{-1}^1 |f g |\;d \mu:f\in X, \Vert f\Vert_X\leq 1\right\}$. 
For the above notions and claims see \cite[Ch.16, \S 112]{ZII}, for example.

\begin{lemma}\label{lemma2.1}
The following conditions for a B.f.s.\ $X$ over $(-1,1)$ with the Fatou property  are equivalent.
\begin{enumerate}[(a)]
\item $X$ has o.c. norm.
\item $X$ is separable.
\item $X'=X^*$.
\item $X$ is weakly sequentially complete.
\item $X$ does not contain an isomorphic copy of $c_0$.
\item $X$ does not contain an isomorphic copy of $\ell^\infty$.
\end{enumerate}
\end{lemma}

\begin{proof}
Recall that $X$ is always Dedekind complete.
\begin{itemize}
\item[(a)$\Leftrightarrow$(b)] Apply \cite[Corollary I.5.6]{BS} together 
with the fact that $\mu$ is separable in the sense of \cite[Definition I.5.4]{BS}.
\item[(a)$\Leftrightarrow$(c)] See \cite[Corollary I.4.3]{BS}.
\item[(c)$\Rightarrow$(d)] See \cite[Corollary I.5.3]{BS}.
\item[(d)$\Leftrightarrow$(e)] See \cite[Theorem 1.c.4]{LZII}.
\item[(e)$\Rightarrow$(f)]     Obvious as $\ell^\infty$ contains an isomorphic copy of $c_0$.
\item[(f)$\Rightarrow$(a)]     See \cite[Proposition 1.a.7]{LZII}.
\end{itemize}
\end{proof}

Given $f\in L^0$, its decreasing rearrangement $f^*:\colon [0,2]\to [0,\infty]$ is the right continuous inverse of its distribution function $ \lambda \mapsto \mu\left(\left\{t\in (-1,1):|f(t)|>\lambda \right\}\right)$, \cite[p.39]{BS}. 
Measurable functions $f$ and $g$ are said
to be equimeasurable if they have the same distribution function.
Let $X$ be a B.f.s.\ over $(-1,1)$ \emph{with} the Fatou property. 
We say that $X$ is rearrangement invariant, \cite[p.59]{BS}, \cite[pp.117-118]{LZII}, briefly r.i., 
if, whenever $g\in L^0$ satisfies $g^*\leq f^*$ for some $f\in X$, we necessarily 
have $g\in X$. In this case, its associate space $X'$ is also r.i., \cite[Proposition II.4.2]{BS}. 
Moreover, $L^\infty \subseteq X\subseteq L^1$ with continuous inclusions, \cite[Corollary II.6.7]{BS}.
Even though $X_a$ may fail to  have the Fatou property, this is the \textit{only} 
reason why it fails to be a r.i.\ space. In fact,  given $f\in X_a$ every function 
$g\in L^0$ which is equimeasurable with $f$ also belongs to $X_a$.

Let us recall the lower and upper Boyd indices of a r.i.\  space $X$ over $(-1,1)$. 
Given $t>0$, define its corresponding dilation operator $E_t\colon X\to X$ 
by $E_t(f)(x):= f(t x)$ if $-1<t x <1$ and $E_t(f)(x):= 0$ otherwise. 
Then $E_t$ is bounded, linear and its operator norm $\Vert E_t\Vert_{\mbox{\scriptsize{op}}}$ 
satisfies $\Vert E_t\Vert_{\mbox{\scriptsize{op}}}\leq \max\{t,1\}$. Define the lower and upper 
Boyd indices by
\begin{align*}
\underline{\alpha}_X:= \sup_{0<t<1}\frac{\log 
\Vert E_{1/t}\Vert_{\mbox{\scriptsize{op}}}}{\log t}\hspace{5mm} 
\mbox{ and } \hspace{5mm}\overline{\alpha}_X:= \inf_{1<t<\infty}\frac{\log \Vert E_{1/t}\Vert_{\mbox{\scriptsize{op}}}}{\log t},
\end{align*}
respectively, \cite[Definition III.5.12]{BS}. Then $0\leq \underline{\alpha}_X\leq\overline{\alpha}_X\leq 1$. Moreover,
\begin{equation}\label{2.1}
\underline{\alpha}_{X'}=1- \overline{\alpha}_X \hspace{5mm}\mbox{ and }\hspace{5mm} 
\overline{\alpha}_{X'}=1- \underline{\alpha}_X,
\end{equation}
\cite[Proposition III.5.13]{BS}.

A typical r.i.\  space is $L^p$ for $1\leq p\leq \infty$ in which case 
$\underline{\alpha}_{L^p}=\overline{\alpha}_{L^p}=1/p$. More 
generally, for the Lorentz spaces $L^{p,q}:= L^{p,q}(-1,1)$ we 
have $\underline{\alpha}_{L^{p,q}}=\overline{\alpha}_{L^{p,q}}=1/p$ 
if $1<p<\infty$ and $1\leq q\leq \infty$ or if $p=q= \infty$, \cite[Definition IV.4.4 and Theorem IV.4.6]{BS}.

The following result, which can be found in \cite[Proposition 2.b.3]{LZII}, will be used in the sequel.

\begin{lemma}\label{lemma2.2}
Let $X$ be a r.i.\  space over $(-1,1)$ such that 
$0<\alpha <\underline{\alpha}_X\leq\overline{\alpha}_X<\beta<1$. 
Then there exist $p,q$ satisfying both $1/\beta <p<q<1/\alpha$ and 
$L^q\subseteq X\subseteq L^p$ with continuous inclusions.
\end{lemma}

Important will be the space weak-$L^1$, denoted by $L^{1,\infty}:= L^{1,\infty}(-1,1)$. 
It is not a Banach space but a quasi-Banach space, \cite[Definition IV.4.1]{BS}.

Define the finite Hilbert transform $T(f)$ of each $f\in L^1$ by the Cauchy principal value integral 
\begin{align*}
(T(f))(t):= \lim_{\varepsilon\to 0+} \frac{1}{\pi} \left(\int_{-1}^{t-\eps}+ \int_{t+\eps}^1\right)\frac{f(x)}{x-t}\; d x,
\end{align*}
which exists for a.e.\  $t\in (-1,1)$. Then $T(f)\in L^0$. 
Moreover, by Kolmogorov's Theorem, \cite[Theorem III.4.9(b)]{BS}, $T$ 
maps $L^1$ into $L^{1,\infty}$ and the resulting linear operator $T\colon L^1 \to L^{1,\infty}$ is continuous.

Now, let $X$ be a r.i.\  space over $(-1,1)$ with non-trivial lower and upper Boyd indices, that is, 
\begin{align*}
0<\lbal_X\leq\ubal_X<1.
\end{align*}
A theorem of Boyd, \cite[Theorem III.5.18]{BS}, ensures that the Hilbert transform 
$H$ is bounded from any r.i.\  space over $\R$, with non-trivial lower and upper 
Boyd indices, into itself. Since we can write $T(f)= -\chi_{(-1,1)}H\left(f\chi_{(-1,1)}\right)$ 
for each $f\in L^1$, which has a natural extension to a function on $\R$, Boyd's theorem 
ensures that $T$ also  maps $X$ into $X$ boundedly. The corresponding bounded linear 
operator is denoted by $T_X\colon X\to X$. Moreover, the finite Hilbert transform 
operator $T_{X'}\colon X'\to X'$ is then also defined; see (\ref{2.1}).

\begin{lemma}\label{lemma2.3}
Let $X$ be a r.i.\  space over $(-1,1)$ satisfying $0<\lbal_X\leq \ubal_X <1$. 
Then $\mbox{sim}\Bb\subseteq X_a$ and hence, $X_a=X_b$. Moreover, 
the inclusion $T_X(X_a)\subseteq X_a$ holds.
\end{lemma}

\begin{proof}
Via Lemma \ref{lemma2.2}, select $1<q<\infty$ satisfying $L^q\subseteq X$ 
continuously. Then $\mbox{sim}\Bb\subseteq L^q=\left(L^{q}\right)_a\subseteq X_a$, 
where the inclusion $(L^q)_a \subseteq X_a$ is a consequence of the continuity of 
the natural injection from $L^q$ into $X$. So, $X_a=X_b$ as discussed after the definition of $X_b$.

Next, given $A\in \Bb$, we have $T_X(\chi_A)=T_{L^q}(\chi_A)\in L^q\subseteq X_a$ and hence,
\begin{align*}
T_X(X_a)=T_X(X_b)=T_X(\overline{\mbox{sim}\Bb})\subseteq \overline{T_X(\mbox{sim}\Bb)}\subseteq X_a.
\end{align*}
\end{proof}

Let $X$ be a r.i.\  space over $(-1,1)$ satisfying $0<\lbal_X\leq \ubal_X <1/2$ or 
$1/2<\lbal_X\leq \ubal_X <1$. Then the operator $T_X\colon X\to X$ is a Fredholm 
operator; see Lemma \ref{lemma2.4} and Lemma \ref{lemma2.5} below. These two 
lemmas are part of Theorems 3.2 and 3.3  in \cite{COR}, respectively. These theorems 
generalize the analogous results when $X:= L^p$ for $1<p<2$ and $2<p<\infty$ 
(recall that $\lbal_{L^p} =\ubal_{L^p}=1/p$), which were originally given in \cite[Theorem 13.9]{J}; 
see \cite[Propositions 2.4 and 2.6]{OE} together with \cite[\S 4-3]{T} for alternative proofs. 
Let $I_X\colon X\to X$ denote the identity operator.

Define the even function $w$ by $w(x):= \sqrt{1-x^2}$ for $x\in (-1,1)$. 
Assume first that $1/2<\lbal_X\leq \ubal_X <1$. Given $f\in X$, define the measurable function 
\begin{equation}\label{2.2}
\widehat{T}_X(f):= - T_X(f w)/w.
\end{equation}

\begin{lemma}\label{lemma2.4}
Let $X$ be a r.i.\  space over $(-1,1)$ satisfying $1/2<\lbal_X\leq \ubal_X <1$.
\begin{enumerate}[(i)]
\item The function $1/w\in X$ and $\ker (T_X)= 	\textnormal{span}\left\{1/w\right\}$.
\item The linear operator $\widehat{T}_X$ defined by (\ref{2.2}) maps $X$ 
boundedly into $X$ and satisfies $T_X\widehat{T}_X=I_X$. In particular, $T_X$ is surjective.
\item The identity $\widehat{T}_X T_X =I_X-P_X$ holds, with $P_X$ 

denoting the bounded linear projection on $X$ given by 
\begin{align*}
f\mapsto P_X(f):= \left(\frac{1}{\pi}\int_{-1}^1 f \;d\mu\right)\frac{1}{w},\hspace{5mm} f\in X.
\end{align*}
\end{enumerate}
\end{lemma}

In the case when $0<\lbal_X\leq \ubal_X <1/2$, define for every $f\in X$ 
a measurable function $\widecheck{T}_X(f)$ by
\begin{equation}\label{2.3}
\widecheck{T}_X(f):= -w T\left(f/w\right).
\end{equation}
Note that $f/w\in L^1$ because $1/w\in X'$.

\begin{lemma}\label{lemma2.5}
Let $X$ be a r.i.\  space over $(-1,1)$ satisfying $0<\lbal_X\leq \ubal_X <1/2$.
\begin{enumerate}[(i)]
\item The linear operator $T_X\colon X\to X$ is injective.
\item The linear operator $\widecheck{T}_X$ defined by (\ref{2.3}) 
maps $X$ boundedly into $X$ and satisfies $\widecheck{T}_X T_X =I_X$. 
\item The range of $T_X$ is the closed linear subspace of $X$ given by
\begin{align*}
\mathscr{R}(T_X)=\left\{h\in X: \int_{-1}^1\frac{h}{w} \;d \mu =0\right\}
\end{align*}
and $(T_X)^{-1}=\widecheck{T} |_{\mathscr{R}(T_X)}$, that is, 
$T_X\left(\widecheck{T}_X(h) \right)=h$ for all $h\in \mathscr{R}(T_X)$.
\end{enumerate}
\end{lemma}

Let $X$ be a r.i.\  space over $(-1,1)$ satisfying $0<\lbal_X\leq \ubal_X <1$. 
In \cite[\S 4]{COR}, a B.f.s.\  $[T,X]$ over $(-1,1)$ with the Fatou property has 
been presented in order to obtain a \emph{maximal} continuous, linear, $X$-valued 
extension of the finite Hilbert transform operator $T_X\colon X\to X$. 
According to \cite[Theorem 4.6]{COR}, the space $[T,X]$ is the \emph{largest} 
B.f.s.\ within $L^1$, containing $X$, to which $T_X$ admits an $X$-valued, 
continuous linear extension. We shall improve this by showing that $[T,X]$ is 
actually the largest one amongst such B.f.s.' within $L^0$, that is, we can omit 
the condition that such B.f.s.' are contained in $L^1$. This will be presented in 
Proposition \ref{proposition2.7} below. To be precise, recall from \cite{COR} that 
\begin{equation}
\label{2.4}[T,X]:= \left\{f\in L^1:T(h)\in X \mbox{ for all }h\in L^0 \mbox{ with } |h|\leq |f| \mbox{ a.e.}\right\}.
\end{equation}
Then $[T,X]$ is an order ideal of $L^1$ and hence, of $L^0$. 
We now recall various properties of $[T,X]\subseteq L^1$ which will be used in the sequel.

\begin{lemma}\label{lemma2.6}
Let $X$ be a r.i.\  space over $(-1,1)$ satisfying $0<\lbal_X\leq \ubal_X <1$. 
Then the following statements hold.
\begin{enumerate}[(i)]
\item The following conditions are equivalent for a function $f\in L^1$.
\begin{enumerate}[(a)]
\item $f\in [T,X]$.
\item $\sup_{|h|\leq |f|}\Vert T(h)\Vert_X<\infty.$
\item $T(f\chi_A)\in X$ for all $A\in \Bb$.
\item $fT_{X'}(g)\in L^1$ for every $g\in X'$.
\end{enumerate}
\item Define
\begin{align*}
\Vert f\Vert_{[T,X]}:= \sup_{|h|\leq |f|}\Vert T(h)\Vert_X, \hspace{5mm} f\in [T,X].
\end{align*}
Then $\Vert \cdot\Vert_{[T,X]}$ is a lattice norm for which $[T,X]$ is a 
B.f.s.\  over $(-1,1)$ with the Fatou property. Moreover, the continuous inclusion $X\subseteq [T,X]$ holds.
\item The following Parseval type formula is valid:
\begin{align*}
\int_{-1}^1 f T_{X'}(g)\; d \mu =-\int_{-1}^1 g T(f) \;d \mu ,\hspace{5mm} f\in [T,X],~ g\in X'.
\end{align*}
\end{enumerate}
\end{lemma}

\begin{proof}
For (i), see \cite[Proposition 4.1]{COR}. Proposition 4.5 of \cite{COR} 
gives the first part of (ii) whereas the continuous inclusion $X\subseteq [T,X]$ is 
immediate by the boundedness of $T_X$. Finally, (iii) is exactly Lemma 4.3 of \cite{COR}.
\end{proof}
\begin{proposition}\label{proposition2.7}
Let $X$ be a r.i.\  space over $(-1,1)$ satisfying $0<\lbal_X\leq \ubal_X <1$. 
Then $[T,X]$ is the largest B.f.s.\ in $L^0$ containing $X$, to 
which $T_X\colon X\to X$ admits an $X$-valued, continuous linear extension.
\end{proposition}
\begin{proof}
Take any B.f.s.\ $Z\subseteq L^0$ for which $X\subseteq Z$ 
and such that there is a continuous linear extension $S\colon Z\to X$ of $T_X$. 
Let $f\in Z$. Select a sequence $\{s_n\}_{n=1}^\infty \subseteq \si \Bb$ satisfying
\begin{equation}\label{2.5}
s_n\geq 0 \mbox{ for } n\in \N\mbox{ and } s_n\uparrow |f| \mbox{ pointwise}.
\end{equation}
For each $n\in \N$ we claim that 
\begin{align}\label{2.6}
\notag\frac{1}{2}\Vert \chi_{(0,1)}\Vert_X \int_{-1}^0 s_n(x)\;d x 
&\leq \Vert T_X(s_n\chi_{(-1,0)})\Vert_X=\Vert S(s_n\chi_{(-1,0)})\Vert_X\\
&\leq \Vert S\Vert_{\mbox{\tiny{op}}} \Vert s_n\Vert_Z\leq \Vert S\Vert_{\mbox{\tiny{op}}} \Vert f\Vert_Z.
\end{align}
In fact, the first inequality can be obtained via the arguments 
in \cite[Proof of Lemma 4.4]{COR} as follows. Given $0<t<1$, we have 
\begin{align*}
|T_X(s_n\chi_{(-1,0)})(t)|=\int_{-1}^0\frac{s_n(x)}{|x-t|}\; d x \geq \frac{1}{2} \int_{-1}^0s_n(x) \;d x
\end{align*}
because $|x-t|\leq 2$ for all $x\in (-1,0)$. Since $\Vert \cdot \Vert_X$ is a lattice norm, this yields
\begin{align*}
\Vert T_X(s_n\chi_{(-1,0)})\Vert_X \geq \Vert T_X(s_n\chi_{(-1,0)})\chi_{(0,1)}
\Vert_X\geq \left(\frac{1}{2} \int_{-1}^0s_n(x)\; d x\right)\Vert \chi_{(0,1)}\Vert_X,
\end{align*}
as claimed. The rest of (\ref{2.6}) is valid by the fact that $T_X(s_n\chi_{(-1,0)})
=S(s_n \chi_{(-1,0)})$ in $X$ together with the definition of the operator 
norm $\Vert S\Vert_{\mbox{\tiny{op}}}$ and (\ref{2.5}).

By the Monotone Convergence Theorem, we have from (\ref{2.6}) that
\begin{align*}
\frac{1}{2}\Vert \chi_{(0,1)}\Vert_X\int_{-1}^0 |f(x)|\; d x = \frac{1}{2}\Vert \chi_{(0,1)}\Vert_X 
 \lim_{n\to\infty} \int_{-1}^0 s_n(x)\; d x\leq \Vert S\Vert_{\mbox{\tiny{op}}} \Vert f\Vert_Z,
\end{align*}
so that $\int_{-1}^0 |f(x)|\; d x<\infty$. By an analogous argument 
we can verify that also $\int_0^1 |f(x)|\; d x<\infty$. So, $f\in L^1$.

Thus, we have established that $Z\subseteq L^1$, which enables 
us to apply \cite[Theorem 4.6]{COR} to conclude that that $Z\subseteq [T,X]$. This completes the proof.
\end{proof}

In view of Proposition \ref{proposition2.7} above, let us call $[T,X]$ 
the optimal lattice domain for the operator $T_X\colon X\to X$ following the terminology of \cite[\S 3]{CR-N}.
Note that $[T,X]\subseteq L^1$ is a \textit{proper} containment. Otherwise the continuous extension
$T\colon[T,X]\to X$ of $T_X$ would imply that $T\colon L^1\to X$ is continuous which is 
not the case as $T(L^1)\not\subseteq L^1$.

We conclude this section with some sufficient conditions for the equality $X=[T,X]$.

\begin{lemma}\label{lemma2.8}
The following statements hold.
\begin{enumerate}[(i)]
\item Let $X$ be a r.i.\  space satisfying $0<\lbal_X\leq \ubal_X <1/2$ or $1/2<\lbal_X\leq \ubal_X <1$. Then $[T,X]=X$.
\item If $1<p<\infty$, then $[T,L^p]=L^p$.
\end{enumerate}
\end{lemma}

\begin{proof}
(i) See Theorem 4.7 of \cite{COR} whose proof explicitly yields $X=[T,X]$.

(ii)  For $1<p<2$ or $2 < p < \infty$, we can apply part (i) as 
$\lbal_{L^p}=\ubal_{L^p}=1/p$. When $p=2$ see Theorem 5.3 in \cite{COR}; its proof explicitly verifies $L^2=[T,L^2]$.
\end{proof}


\section{Vector measure arising from the finite Hilbert transform}\label{section3}

We shall proceed with the notation from Section  \ref{section2}. 
Throughout this section let $X$ denote a r.i.\  space over $(-1,1)$ 
with non-trivial Boyd indices and let $T_X\colon X\to X$ denote the 
corresponding finite Hilbert transform operator. The $X$-valued set function
\begin{align*}
m_X\colon A\mapsto T_X(\chi_A),\hspace{5mm} A\in \Bb,
\end{align*}
turns out to be a $\sa$-additive vector measure. This provides the means 
to investigate aspects of $T_X$ from a different point of view, including its 
extensions to spaces other than $[T,X]$. Our main result is Theorem 
\ref{theorem3.10} regarding the relationships between $X,[T,X]$ and 
the spaces $L^1(m_X)$ and $L^1_w(m_X)$ of $m_X$-integrable, 
respectively scalarly $m_X$-integrable functions.

Let us begin with an arbitrary finitely additive set function $\nu\colon \Bb\to X$; its 
range in $X$ is denoted by by $\Rr(\nu):=\{\nu(A):A\in \Bb\}$. By $|\nu| $ we 
denote its variation as defined in the scalar case, \cite[Definition I.1.4]{DU}. 
Then $|\nu|\colon \Bb\to [0,\infty]$ is an extended real-valued, finitely additive 
set function. If $\nu$ happens to be $\sa$-additive, then so is $|\nu|$ 
(see Proposition I.1.9 and its proof in \cite{DU}). Given $x^*\in X^*$, 
let $\langle \nu,x^*\rangle \colon \Bb \to \C$ denote the finitely additive 
set function $A\mapsto \langle\nu(A),x^*\rangle$ for $A\in \Bb$. 
Define its semivariation $\Vert\nu \Vert \colon \Bb\to [0,\infty]$ by 
\begin{align*}
\Vert \nu\Vert (A):= \sup\left\{|\langle \nu,x^*\rangle|(A):x^* \in X^*,\Vert x^*\Vert_{X^*}\leq 1\right\}, \hspace{5mm} A\in \Bb,
\end{align*}
\cite[p.2]{DU}. It follows that $\Vert \nu\Vert (\Bb)\subseteq [0,\infty)$ if and 
only if $\Rr(\nu)$ is a bounded set in $X$. This fact is part of Proposition I.1.11 in \cite{DU}, as is the formula 
\begin{equation}\label{3.1}
\Vert \nu\Vert (A) =\sup \left\Vert \sum_{j=1}^n a_j \nu(A_j\cap A)\right\Vert_X,\hspace{5mm} A\in \Bb,
\end{equation}
where the supremum is taken over all pairwise disjoint sets $A_1,\ldots ,A_n\in \Bb$ 
and scalars $a_1,\ldots, a_n\in \C$ satisfying $|a_j|\leq 1$ for $j=1,\ldots,n$, with $n\in \N$.

Recall the B.f.s.\ $[T,X]\subseteq L^1$; see (\ref{2.4}) together with its 
properties from Lemma \ref{lemma2.6}. Given $f\in [T,X]$, define a finitely 
additive set function $\nu_f\colon \Bb\to X$ by
\begin{equation}\label{3.2}
\nu_f(A):= T(f\chi_A),\hspace{5mm} A\in \Bb;
\end{equation}
see Lemma \ref{lemma2.6}(i). Then $\Rr(\nu_f)$ is bounded in $X$ 
as $\sup_{A\in \Bb}\Vert \nu_f(A)\Vert_X\leq \Vert f \Vert_{[T,X]}<\infty$ via Lemma \ref{lemma2.6}(ii).

\begin{lemma}\label{lemma3.1}
Let $X$ be a r.i.\  space over $(-1,1)$ satisfying $0<\lbal_X\leq \ubal_X <1$ 
and $f\in [T,X]$. Then $\Vert f\chi_A\Vert_{[T,X]} =\Vert \nu_f\Vert (A), A\in \Bb$.
\end{lemma}

\begin{proof}
Fix $A\in \Bb$. Then the linear operator $\Phi_{f,A}\colon L^\infty\to X$ 
given by $\Phi_{f,A}(\phi):= T(\phi f \chi_A)$ is bounded because, for each non-zero $\phi\in L^\infty$, we have
\begin{align*}
\Vert \Phi_{f,A}(\phi)\Vert_X=\Vert T(\phi f \chi_A)\Vert_X=\Vert \phi\Vert_\infty 
\Vert T(f \chi_A \phi/\Vert \phi\Vert_\infty)\Vert_X\leq \Vert f\Vert_{[T,X]} \Vert \phi\Vert_\infty;
\end{align*}
see Lemma \ref{lemma2.6}(ii). Moreover, $\Vert f \chi_A\Vert_{[T,X]}= \Vert \Phi_{f,A}\Vert_{\mbox{\tiny{op}}}$ because
\begin{align*}
\Vert f\chi_A \Vert_{[T,X]}:= \sup_{|h|\leq |f\chi_A|}\Vert T(h)\Vert_X=\sup_{|\phi|\leq 1}\Vert T(\phi f \chi_A)\Vert_X.
\end{align*}
On the other hand, by (\ref{3.1}) with $\nu_f$ in place of $\nu$, we have
\begin{equation}\label{3.3}
\Vert\nu_f\Vert(A) =\sup\{ \Vert T(s f\chi_A)\Vert_X: s\in \si\Bb,|s|\leq 1\}.
\end{equation}
Since the subset $\{s\in \si\Bb,|s|\leq 1\}$ is dense in the closed unit ball 
of $L^\infty$, the right side of (\ref{3.3}) equals the operator norm of $\Phi_{f,A}$. This establishes the lemma.
\end{proof}

The $X$-valued set function $m_X\colon A\mapsto T_X(\chi_A)$ on $\Bb$ is 
clearly finitely additive. We can also write $m_X=\nu_{\bf{1}}$ in the notation 
(\ref{3.2}) with $f$ the function $\bf{1}$ which is constantly equal to one. 
Let us see that $m_X$ is a vector measure, that is, it is $\sa$-additive. 
Lemma \ref{lemma2.3} shows that $\chi_A\in X_a$ for all $A\in \Bb$. 
So, if $A_n\downarrow\emptyset$ in $\Bb$, then $\chi_{A_n}\downarrow 0$ 
in $X_a$ and so order continuity of the norm in $X_a$ implies that 
$\lim_{n\to \infty}\chi_{A_n} =0$ in $X_a$ (hence, also in $X$). 
The boundedness of $T_X$ then ensures that $\lim_{n\to \infty} m_X(A_n)
=\lim_{n\to \infty} T_X( \chi_{A_n})=0$. Thus, $m_X$ is $\sa$-additive.

We shall say that $m_X$ is the vector measure \emph{arising from} the 
operator $T_X\colon X\to X$. In the terminology of \cite[Proposition 4.4]{ORS}, 
$m_X$ is the vector measure associated with the restriction 
$T_X|_{X_a}\colon X_a\to X$ of $T_X$ to the B.f.s.\ $X_a$  with o.c.\  norm. 
Note that $X_a$ may not have the Fatou property; see Remark \ref{remark3.12}(b-2) below. 
We shall again discuss this and its implications after Proposition \ref{proposition3.2} below.

A set $A\in \Bb$ is said to be $m_X$-null if $|\langle m_X,x^*\rangle |(A)=0$ 
for every $x^*\in X^*$. Note that $A\in \Bb$ is $m_X$-null if and only if 
$m_X(B)=0$ for all $B\in \Bb$ with $B\subseteq A$. A property which 
holds outside an $m_X$-null set is said to hold $m_X$-almost 
everywhere, briefly $m_X$-a.e., as in the case of scalar measures.

Let us present some relevant properties of $m_X$ in the following result. 
Properties (v)-(vii), which indicate that $m_X$ is a rather non-trivial vector 
measure, are due to the fact that the kernel generating $T_X$ is far from 
being non-negative and has no monotonicity properties.

\begin{proposition}\label{proposition3.2}
Let $X$ be a r.i.\  space over $(-1,1)$ satisfying $0<\lbal_X\leq \ubal_X <1$. 
The following statements hold for the vector measure $m_X\colon \Bb\to X$ 
arising from the operator $T_X\colon X\to X$.
\begin{enumerate}[(i)]
\item The range $\Rr(m_X)\subseteq X_a$.
\item For every $g\in X'\subseteq X^*$, the scalar measure $\langle m_X, g\rangle$ is given by 
\begin{align*}
\langle m_X, g\rangle(A)=-\int_A T_{X'} (g)\, d \mu,\hspace{5mm} A\in \Bb.
\end{align*}
\item There exists $g_0\in X'$ satisfying $|\langle m_X, g_0\rangle|=\mu$ on $\Bb$.
\item A set $A\in\Bb$ is $\mu$-null if and only if it is $m_X$-null.
\item The vector measure $m_X$ has totally infinite variation. 
In other words, $|m_X|(A)=\infty$ for every non-$m_X$-null set $A\in \Bb$.
\item $\Rr(m_X)$ is not a relatively compact subset of $X$.
\item $\Rr(m_X)$ is not an order  bounded subset of $X$.
\end{enumerate}
\end{proposition}

\begin{proof}
(i) This follows from Lemma \ref{lemma2.3} as $\Rr(m_X)\sub T_X(\si\Bb)\sub T_X(X_a)\sub X_a$.

(ii) Apply Lemma \ref{lemma2.6}(iii) with $f:=\chi_A$ for each $A\in \Bb$ to obtain
\begin{align*}
\langle m_X, g\rangle(A)&=\langle m_X(A), g\rangle=\int_{-1}^1 
T_{X}(\chi_A) g\, d \mu=-\int_{-1}^1 \chi_A T_{X'}(g) \, d \mu\\
&=-\int_A  T_{X'}(g) \, d \mu.
\end{align*}

(iii) Via Lemma \ref{lemma2.2} there exists $1<p<2$ such that 
$X\sub L^p$ continuously. Since the conjugate index $p':= p/(p-1)$ 
satisfies $2<p'<\infty$, we can apply Lemma \ref{lemma2.5}(iii) with 
$X:= L^{p'}$ and the fact that $w$ is an even function to see that the 
sign function $\sa:= -\chi_{(-1,0)}+\chi_{(0,1)}$ belongs to $\Rr(T_{L^{p'}})$ 
and that $T_{L^{p'}}(g_0) =\sa$, where $g_0\in L^{p'}$ is given by
\begin{equation}\label{3.4}
g_0:= \widecheck{T}_{L^{p'}}(\sa) 
=-w T \left(\sa/w\right).
\end{equation}
Since $L^{p'}=(L^p)'\sub X'$, it follows that $g_0\in X'$ satisfies $T_{X'}(g_0) =\sa$. 
Then $\langle m_X,g_0\rangle (A) =-\int_A T_{X'}(g_0)\, d \mu = -\int_A \sa \;d \mu$ 
for $A\in \Bb$ via part (ii) with $g:= g_0$, by which 
$|\langle m_X,g_0\rangle| (A)=\int_A |\sa| \;d \mu=\mu(A)$ for $A\in \Bb$. So, (iii) is verified.

(iv) Assume that $A\in \Bb$ is $\mu$-null. Then every $B\in \Bb$ with $B \sub A$ 
satisfies $\mu(B)=0$ and hence, $\chi_B=0$ a.e.\  so that $m_X(B)=T_X(\chi_B)=0$. Thus, $A$ is $m_X$-null.

Conversely, assume that $A$ is $m_X$-null, in which case it is also 
$|\langle m_X,g \rangle|$-null for all $g\in X'$. In particular, $A$ is then 
$\mu$-null as $\mu(A)=|\langle m_X,g_0 \rangle|(A)=0$ with $g_0\in X'$ as in part (iii).

(v) Select $1<p<2$ as in the proof of part (iii), in which case $1/w\in L^p$. 
Consider the $L^p$-valued vector measure $\eta\colon A\mapsto \chi_A$ on 
$\Bb$. Then there exist constants $C_1,C_2>0$ such that
\begin{equation}\label{3.5}
|\eta|(A) \leq C_1 |m_X|(A) + C_2\mu(A), \hspace{5mm} A\in \Bb.
\end{equation}
Indeed, let $j_p\colon X\to L^p$ denote the natural embedding. 
As $m_{L^p}= j_p\circ m_X$ we have, from Lemma \ref{lemma2.4}(iii) with $X:= L^p$, that 
\begin{align}\label{3.6}
\notag\chi_A&=\left(\widehat{T}_{L^p} T_{L^p}\right)(\chi_A)+ P_{L^p}(\chi_A)
= \left(\widehat{T}_{L^p}\circ m_{L^p}\right)(A)+P_{L^p}(\chi_A)\\
&= \left(\widehat{T}_{L^p}\circ j_p\circ m_X\right)(A) +P_{L^p}(\chi_A),\hspace{5mm}A\in \Bb.
\end{align}
Setting $C_1:= \Vert \widehat{T}_{L^p}\circ j_p\Vert_{\mbox{\tiny{op}}}$ 
and $C_2:= \Vert 1/(\pi w)\Vert_{L^p}$, it follows from (\ref{3.6}) and the identity 
$P_{L^p}(\chi_A)=(\mu(A)/\pi)\cdot (1/w)$ that
\begin{align*}
\Vert \eta (A)\Vert_{L^p}\leq C_1 \Vert m_X(A)\Vert_X+C_2\mu(A),\hspace{5mm} A\in \Bb.
\end{align*}
This verifies (\ref{3.5}). 

Now, given a non-$m_X$-null set $A\in\Bb$, we have from part (iv) 
that $A$ is not $\mu$-null. So, $|\eta|(A)=\infty$, \cite[Example I.1.16]{DU}. This and (\ref{3.5}) establish part (v).

(vi) We use the notation from the proof of part (v). Let $\eta_1$ 
denote the $L^p$-valued vector measure $A\mapsto P_{L^p}(\chi_A)
=\mu(A) \cdot1/(\pi w)$ on $\Bb$. Then we can rewrite (\ref{3.6}) as 
\begin{align*}
\eta(A)= \left(\widehat{T}_{L^p}\circ j_p\circ m_X\right)(A) +\eta_1(A),\hspace{5mm} A\in \Bb.
\end{align*}
Assume that $\Rr(m_X)$ is a relatively compact subset of $X$. Then the 
vector measure $\widehat{T}_{L^p}\circ j_p\circ m_X\colon \Bb\to L^p$ has 
relatively compact range in $L^p$ and hence, so does $\eta$ because $\Rr(\eta_1)$ 
is clearly compact (as $\eta_1$ takes its values in a 1-dimensional subspace of $L^p$). 
But, $L^p\sub L^1$ continuously and so it follows from \cite[Example III.1.2]{DU} 
that $\eta$ does \emph{not} have relatively compact range. This contradiction verifies part (vi).

(vii) Assume, on the contrary, that $m_X$ does have order bounded range, 
that is, there is a positive function $h$ in $X$ such that 
\begin{equation}\label{3.7}
|m_X(A)|\leq h\; \mbox{ a.e.}, \hspace{5mm} A\in \Bb.
\end{equation}
Apply Lusin's Theorem to $h\in L^0$ to find a compact subset $K\sub(-1,1)$ 
such that $\mu((-1,1)\setminus K)<1$ and the restriction $h|_K$ of $h$ to $K$ 
is a continuous function. Then $M:= \max_{t\in K} |h(t)|<\infty$.

Given $t\in K$ we have
\begin{align*}
\lim_{x\to t} \left|m_X((t,1))(x)\right|=\lim_{x\to t} \left|T_X(\chi_{(t,1)})(x)\right|
=\frac{1}{\pi} \lim_{x\to t} \left|\ln \left|\frac{1-x}{t-x}\right|\,\right|=\infty,
\end{align*}
because $T_X(\chi_{(t,1)})(x)=(1/\pi)\ln \left|\tfrac{1-x}{t-x}\right|$ for $x\in (-1,1)\setminus\{t\}$ 
which can be directly verified via the definition of the principal value integral. 
So, there exists an open subinterval $U(t)\sub (-1,1)$ containing $t$ such that 
\begin{equation}\label{3.8}
\left|m_X((t,1))(x)\right|>2M,\hspace{5mm} x\in U(t)\setminus\{t\}.
\end{equation}
Now, since $\{U(t)\}_{t\in K}$ is an open cover of $K$, we can find finitely many 
points $t_1,\ldots, t_n\in K$ such that  $K\sub \bigcup_{j=1}^n U(t_j)$. 
Observe that $\mu(K)>0$ as $\mu((-1,1)\setminus K)<1$. 
Accordingly, $\mu(U(t_j)\cap K)>0$ for some $j\in \{1,\ldots,n\}$. Setting $t:= t_j$ in (\ref{3.8}) yields 
\begin{equation}\label{3.9}
\left|m_X((t_j,1))(x)\right|>2M, \hspace{5mm} x\in \left( U(t_j)\cap K\right)\setminus\{t_j\}.
\end{equation}
However, (\ref{3.9}) contradicts (\ref{3.7}) because (\ref{3.7}), for $A:= \left( U(t_j)\cap K\right)\setminus\{t_j\}$, gives
\begin{align*}
\left|m_X((t_j,1))(x)\right|\leq h(x)\leq M, \hspace{5mm}\mbox{ a.e. } x\in \left( U(t_j)\cap K\right)\setminus\{t_j\},
\end{align*}
with $\left( U(t_j)\cap K\right)\setminus\{t_j\}$ having positive measure. 
Therefore, $\Rr(m_X)$ fails to be order bounded.
\end{proof}

\begin{remark}\label{remark3.3}
(i) The function $g_0\in X'\sub X^*$ in (\ref{3.4}) is called a \emph{Rybakov functional} 
for $m_X$ (cf. \cite[Theorem IX.2.2]{DU}) because $m_X$ and $|\langle m_X,g_0\rangle |$ 
have the same null sets (by (ii) and (iii) of Proposition \ref{proposition3.2}).

(ii) Proposition \ref{proposition3.2}(v) implies that $L^1(|m_X|)=\{0\}$ and hence, 
$L^1(|m_X|)\neq L^1(m_X)$. Since $L^1(|\langle m_X,g_0\rangle |)=L^1$ 
(see Proposition \ref{proposition3.2}(iii)), it follows from Corollary 3.19(ii) 
of \cite{ORS} that the inclusion $L^1(m_X)\sub L^1$ is always \emph{proper}.

(iii) By part (i) we see that $m_X$ is absolutely continuous with respect to 
$\mu$, \cite[Theorem I.2.1]{DU}, and so the question arises of whether 
$m_X$ possesses an $X$-valued Pettis integrable density with respect 
to $\mu$. The answer is negative. If such a density existed, then $m_X$ 
would have $\sa$-finite variation, \cite[Proposition 5.6(iv)]{vD}, which 
is impossible by Proposition \ref{proposition3.2}(v) above and the fact that $\mu$ is non-atomic.
\end{remark}

Recall that a bounded linear operator between Banach spaces is 
\emph{completely continuous} if it maps relatively weakly compact 
sets to relatively compact sets. Clearly every compact operator is completely 
continuous. A linear operator between Banach lattices is called 
\emph{order bounded} if it maps order bounded sets to order bounded sets.

\begin{corollary}\label{corollary3.4}
Let $X$ be a r.i.\  space over $(-1,1)$ satisfying $0<\lbal_X\leq \ubal_X <1$. 
The operator $T_X\colon X\to X$ is neither completely continuous nor order 
bounded. In particular, $T_X$ is not a compact operator.
\end{corollary}

\begin{proof}
Let us first show that $T_X$ is not completely continuous. The $X$-valued 
set function $\eta\colon  A\mapsto \chi_A$ on $\Bb$ is $\sa$-additive 
because $\chi_A\in X_a$ for all $A\in\Bb$ and $X_a$ has o.c. norm 
(see Lemma \ref{lemma2.3}). Thus, $\Rr(\eta)$ is a relatively  
weakly compact set in $X$, \cite[Corollary I.2.7]{DU}.  On the 
other hand, $T_X(\Rr(\eta))=\Rr(m_X)$ is \emph{not} relatively 
compact by Proposition \ref{proposition3.2}(vi). Therefore, $T_X$ is not completely continuous.

Since $T_X$ maps the order bounded set $\{\chi_A:A\in\Bb\}$ in $X$ to the 
non-order bounded set $\{T(\chi_A):A\in\Bb\}=\Rr(m_X)$ in $X$ 
(see Proposition \ref{proposition3.2}(vii)), the operator $T_X$ is not order bounded.
\end{proof}

\begin{remark}\label{remark3.5}
Let $X$ be as in Corollary \ref{corollary3.4}. 
According to Lemma \ref{lemma2.6}(ii) the natural 
inclusion $\kappa_X\colon X\to [T,X]$ is continuous. 
Moreover, the linear map $\widetilde{T}_X\colon [T,X] \to X$ 
defined by $\widetilde{T}_X(f):= T(f)$ for $f\in[T,X]$, is a continuous, 
$X$-valued extension of $T_X\colon X\to X$; see the discussion 
prior to Lemma \ref{lemma2.6} and part (ii) of Lemma \ref{lemma2.6}. Then the factorization 
\begin{align*}
T_X=\widetilde{T}_X\circ \kappa_X
\end{align*}
together with Corollary \ref{corollary3.4} implies 
that neither $\kappa_X$ nor $\widetilde{T}_X$ is a compact operator.
\end{remark}

A function $f\in L^0$ is called $m_X$-integrable if it satisfies the following two conditions:
\begin{enumerate}
\item[(I-1)] $f$ is integrable with respect to the 
complex measure $\langle m_X, x^*\rangle $ for every $x^*\in X^*$; and
\item[(I-2)] given $A\in \Bb$, there is a unique element $\int_A f\;d m_X\in X$ satisfying
\begin{align*}
\left\langle \int_A f\;d m_X,x^*\right\rangle =\int_A f\;
d\langle  m_X,x^*\rangle ,\hspace{5mm} x^*\in X^*.
\end{align*}
\end{enumerate}

Let $L^1(m_X)$ denote the vector space of all $m_X$-integrable 
functions. We shall identify those $m_X$-integrable  functions which 
are equal $m_X$-a.e. Then the functional
\begin{align*}
f\mapsto \Vert f\Vert_{m_X}:=\sup\left\{\int_{-1}^1 |f|\; 
d |\langle m_X,x^*\rangle | : x^*\in X^*, \Vert x^*\Vert_{X^*}\leq 1\right\}
\end{align*}
on $L^1(m_X)$ defines a norm for which $L^1(m_X)$ is complete and 
$\si \Bb$ is dense in $L^1(m_X)$; see, \cite[Theorems 3.5 and 3.7]{ORS}, 
for example. Given $f\in L^1(m_X)$, its $X$-valued indefinite integral 
$(m_X)_f\colon A\mapsto \int_Af \;d m_X$ on $\Bb$ is $\sa$-additive 
by the Orlicz-Pettis Theorem, \cite[Corollary I.4.4]{DU}. Moreover, the identity 
\begin{equation} \label{3.10}
\Vert (m_X)_f\Vert(-1,1)=\Vert f \Vert_{m_X},
\end{equation}
is valid, \cite[(3.8) on p.106]{ORS}.

Recall from Proposition \ref{proposition3.2}(iv) that the $m_X$-null 
and the $\mu$-null sets coincide. It turns out that $L^1(m_X)$ is an order 
ideal of $L^0$ and that its norm $\Vert \cdot \Vert_{m_X}$ is an o.c. 
lattice norm for which $L^1(m_X)$ is a B.f.s.\ over $(-1,1)$; see, 
for example, \cite[Theorem 3.7]{ORS}. Since the B.f.s.\ $L^1(m_X)$ need 
not have the Fatou property in general (cf. Lemma \ref{lemma3.9} below), 
we cannot apply Lemma \ref{lemma2.1} (with $L^1(m_X)$ in place of $X$) 
to conclude that $L^1(m_X)$ is separable. However, since $\Bb$ is countably 
generated, $m_X$ is a separable vector measure and so the separability of 
$L^1(m_X)$ follows from \cite[Proposition 2]{Ri}. Note that $X$ itself need \emph{not} be separable.

The integration operator $I_{m_X}\colon L^1(m_X)\to X$ associated with the vector measure $m_X$ is defined by
\begin{align*}
I_{m_X}(f):= \int_{-1}^1 f\; d m_X,\hspace{5mm} f\in L^1(m_X).
\end{align*}
Then $I_{m_X}$ is linear, bounded, \cite[p.152]{ORS}, and moreover, satisfies
\begin{equation}\label{3.11}
\Rr(I_{m_X})\sub \overline{\mbox{span}}\,\Rr(m_X)\sub X_a,
\end{equation}
where $\overline{\mbox{span}}\,\Rr(m_X)$ denotes the closed linear span 
of $\Rr(m_X)$ in $X$. In fact, the first inclusion is a consequence of 
$I_{m_X}(\si\Bb)=\mbox{span}\,\Rr(m_X)$ and the denseness of $\si\Bb$ 
in $L^1(m_X)$, whereas the second inclusion follows from Proposition 
\ref{proposition3.2}(i) and the fact that $X_a$ is closed in $X$.

Next, let us establish the continuous inclusion $X_a\sub L^1(m_X)$ via the optimal 
domain theory presented in \cite[\S 4]{ORS}; see also \cite{CR-N}. Recall that 
$T|_{X_{a}}\colon X_a\to X$ denotes the restriction of $T_X$ to $X_a\sub X$. It 
was noted in Section \ref{section2} that $X_a$ is a B.f.s.\ over $(-1,1)$ with o.c.\  norm 
which contains $\si\Bb$ (see Lemma \ref{lemma2.3}). The $X$-valued set function 
$A\mapsto T|_{X_{a}}(\chi_A)$ on $\Bb$ is called the vector measure associated 
with $T|_{X_{a}}$ in the terminology of Proposition 4.4 in \cite{ORS} and is, of 
course, equal to $m_X$. Since the $m_X$-null and the $\mu$-null sets coincide, 
$T|_{X_{a}}$ is $\mu$-determined as formally defined in \cite[p.187]{ORS}. 
This enables us to apply Theorem 4.14 of \cite{ORS} to $T|_{X_{a}}$ to obtain the following fact.

\begin{lemma}\label{lemma3.6}
Let $X$ be a r.i.\  space over $(-1,1)$ satisfying $0<\lbal_X\leq \ubal_X <1$. 
Then $L^1(m_X)$ is the largest amongst all B.f.s.' having o.c. norm into which 
$X_a$ is continuously embedded and to which $T|_{X_{a}}$ admits an 
$X$-valued, continuous linear extension. Further, such an extension is 
unique and equals the integration operator $I_{m_X}\colon L^1(m_X)\to X$. 
In particular, $X_a\sub L^1(m_X)$ continuously and $I_{m_X}(f)= T|_{X_{a}}(f)$ for $f\in X_a$.
\end{lemma}

Let us call $L^1(m_X)$ the o.c.\ optimal lattice domain for the operator 
$T|_{X_{a}}\colon X_a\to X$. It is to be compared with  the optimal lattice 
domain $[T,X]$ of $T_X$ (see Section \ref{section2}).

\begin{remark}\label{remark3.7}
Since $\{\chi_A:A\in \Bb\}$ is a bounded set for the norm $\Vert\cdot\Vert_{\infty}$ 
and $X_a=X_b=\overline{\si \Bb}$ (closure in $X_a$) via Lemma \ref{lemma2.3}, 
it follows that $\{\chi_A:A\in \Bb\}$ is also a bounded set in $X_a$. 
Let $T_{X_a}^{(a)}\colon X_a\to X_a$ be the operator $T|_{X_a}$ interpreted 
as being $X_a$-valued (which is well defined by Lemma \ref{lemma2.3})
and let $m_{X_a}\colon A\mapsto T_{X_a}^{(a)}(A)$ denote its associated vector measure. 
If $T_{X_a}^{(a)}$ was compact, then $\Rr(m_{X_a})=\{T_{X_a}^{(a)}(\chi_A):A\in \Bb\}$ 
would be relatively compact in $X_a$ and hence, also relatively compact in $X$. 
Since this contradicts Proposition \ref{proposition3.2}(vi), we can conclude 
that $T_{X_a}^{(a)}$ is \emph{not} a compact operator. Via Lemma \ref{lemma3.6}, 
the natural inclusion $j_{X_a}\colon X_a\to L^1(m_{X_a})=L^1(m_X)$ 
is continuous and the integration operator $I_{m_X}\colon  L^1(m_X)\to X$ 
is a continuous, $X$-valued extension of $T\mid_{X_a}$. Since $I_{m_X}$ 
takes its values in $X_a$ (see (\ref{3.11})), it follows from the 
factorization $T\mid_{X_a}= I_{m_X}\circ j_{X_a}$ that neither $I_{m_X}$ nor $J_{X_a}$ is compact.
\end{remark}

The B.f.s.\ $L^1(m_X)$ may fail to have the Fatou property. However, 
there always exists a larger B.f.s.\ with the Fatou property which 
contains $L^1(m_X)$ continuously. Such a B.f.s.\ can be realized 
as the space of all scalarly $m_X$-integrable functions. Functions 
of this kind were originally defined in \cite[Definition 2.5]{L1} and a 
systematic study of them was undertaken in \cite{St}.

Let us provide some relevant notions. A function $f\in L^0$ is said to 
be scalarly $m_X$-integrable if it satisfies condition (I-1) above. 
In the special case when $X$ does not have an isomorphic copy of 
$c_0$, every scalarly $m_X$-integrable function is necessarily $m_X$-integrable, 
\cite[Theorem 5.1 and p.302]{L2}.  We shall identify those scalarly 
$m_X$-integrable functions which are equal $m_X$-a.e.\  and denote 
by $L^1_w(m_X)$ the space of (equivalence classes of) all 
scalarly $m_X$-integrable functions. Then $L^1_w(m_X)$ is a Banach space with respect to the norm
\begin{equation}\label{3.12}
\Vert f \Vert_{m_X,w}:= \sup\left\{\int_{-1}^1 |f| \;d |\langle m_X, x^*\rangle|:x^*\in X^*,\Vert x^*\Vert_{X^*}\leq 1\right\}, \hspace{5mm}
\end{equation}
for $f\in L^1_w(m_X)$; see \cite[Theorem 9]{St} for the case of $\R$-valued 
functions and \cite[p.138]{ORS} for $\C$-valued functions. It is clear 
from (\ref{3.12}) and the definition of $\Vert \cdot \Vert_{m_X}$ on 
$L^1(m_X)$ that $L^1(m_X)$ is a closed linear subspace of $L^1_w(m_X)$ 
and that its norm equals that induced by $\Vert \cdot \Vert_{m_X,w}$.

Since the $m_X$-null und $\mu$-null sets coincide, we see that $L^1_w(m_X)$ 
is an order ideal of $L^0$ and that $\Vert \cdot \Vert_{m_X,w}$ is a lattice norm 
for which $L^1_w(m_X)$ is a B.f.s.\ over $(-1,1)$.

\begin{remark}\label{remark3.8}
It is clear that $f\in L^0$ belongs to $L^1_w(m_X)$ if and only if the right-side of 
(\ref{3.12}) is finite. As explained in \cite[p.189]{CR}, this enables us to realize 
$L^1_w(m_X)$ as the B.f.s.\ over $(-1,1)$ relative to $\mu$ corresponding to the 
$[0,\infty]$-valued function norm on $L^0$ given by
\begin{align*}
f\mapsto \sup \left\{\int_{-1}^1 |f| \;d |\langle m_X, x^*\rangle|:x^*\in X^*,\Vert x^*\Vert_{X^*}\leq 1\right\};
\end{align*}
see \cite[Ch. 15]{Z} for function norms and their corresponding B.f.s.'.
\end{remark}

\begin{lemma}\label{lemma3.9}
Let $X$ be a r.i.\  space over $(-1,1)$ satisfying $0<\lbal_X\leq \ubal_X <1$. The following statements are valid.
\begin{enumerate}[(i)]
\item The B.f.s.\ $L^1_w(m_X)$ has the Fatou property and coincides with the bi-associate 
space $L^1(m_X)'':=(L^1(m_X)')'$.
\item The following conditions are equivalent.
\begin{enumerate}[(a)]
\item $L^1(m_X)$ has the Fatou property.
\item $L^1_w(m_X)$ has o.c. norm.
\item $L^1(m_X)=L^1_w(m_X)$.
\item $\si \Bb$ is dense in $L_w^1(m_X)$.
\end{enumerate}
\item The B.f.s.\ $L^1_w(m_X)$ is the minimal B.f.s.\ over $(-1,1)$ with the 
Fatou property which contains (with norm $\leq 1$) the B.f.s.\ $L^1(m_X)$.
\end{enumerate}
\end{lemma}

\begin{proof}
\begin{enumerate}[(i)]
\item See Propositions 2.1 and 2.4 of \cite{CR}.
\item \
\begin{enumerate}[(a)]
\item[(a)$\Leftrightarrow$(b)]See \cite[Proposition 2.3(i)]{CR}. 
\item[(c)$\Rightarrow$(a)]Clear as $L^1_w(m_X)$ always has the Fatou property by part (i).
\item[(b)$\Rightarrow$(c)]Condition (b) gives $(L^1_w(m_X))_a
=L^1_w(m_X)$. Since always $(L^1_w(m_X))_a=L^1(m_X)$, 
\cite[(3.86) on p.145]{ORS}, it is clear that (c) follows.
\item[(c)$\Leftrightarrow$(d)] See \cite[Proposition 3.38(I)]{ORS} with $p:=1$ and $\nu:= m_X$.
\end{enumerate}
\item See \cite[p.191]{CR}.
\end{enumerate}
\end{proof}

Let us present the main result.

\begin{theorem}\label{theorem3.10}
Let $X$ be a r.i.\  space over $(-1,1)$ satisfying $0<\lbal_X\leq \ubal_X <1$. 
The following statements hold for the vector measure 
$m_X\colon A\mapsto T_X(\chi_A)\in X$ on $\Bb$ arising from $T_X$.
\begin{enumerate}[(i)]
\item The natural inclusions
\begin{equation}\label{3.13}
X_a\sub L^1(m_X)\sub [T,X]\sub L^1
\end{equation}
hold and are continuous. Furthermore, the equality
\begin{equation}\label{3.14}
\Vert f\Vert_{m_X}= \Vert f\Vert_{[T,X]},\hspace{5mm} f \in L^1(m_X),
\end{equation}
is valid. In particular, $L^1(m_X)$ is a closed linear subspace of $[T,X]$. Moreover,
\begin{equation}\label{3.15}
I_{m_X}(f \chi_A)=T( f \chi_A),\hspace{5mm} f\in L^1(m_X),\, A\in \Bb.
\end{equation}
\item  The following conditions for a function $f\in [T,X]$ are equivalent.
\begin{enumerate}[(a)]
\item $f\in L^1(m_X)$.
\item $f\in [T,X]_a$.
\item $T(f \chi_A)\in X_a$ for every $A\in \Bb$.
\item The set function $\nu_f\colon \Bb\to X$ given by (\ref{3.2}) is $\sa$-additive.
\end{enumerate}
\item The two spaces $L^1_w(m_X)$ and $[T,X]$ are identical as B.f.s.'.
\end{enumerate}
\end{theorem}

\begin{proof}
(i) For the continuous inclusion $X_a\sub L^1(m_X)$ see Lemma \ref{lemma3.6} above.

To establish the continuous inclusion $L^1(m_X)\sub [T,X]$ 
as well as the identities (\ref{3.14}) and (\ref{3.15}), let us verify the identity 
\begin{equation}\label{3.16}
\Vert s \Vert_{m_X}= \Vert s\Vert _{[T,X]}, \hspace{5mm} s\in \si\Bb.
\end{equation}
From (\ref{3.2}), with $f:= s$, we have 
\begin{equation}\label{3.17}
\nu_s(A)= T ( s \chi_A)= \int_A s \, d m_X, \hspace{5mm} A \in \Bb,
\end{equation}
that is, $\nu_s$ is precisely the indefinite integral $(m_X)_s$ of $s$ 
with respect to $m_X$. Hence, (\ref{3.10}) yields 
$\Vert \nu_s\Vert (-1,1) = \Vert s\Vert_{m_X}$. 
Since $\Vert s \Vert_{[T,X]} = \Vert \nu_s\Vert (-1,1)$ by Lemma \ref{lemma3.1}, we then obtain (\ref{3.16}).

Now let $f\in L^1(m_X)$. Select $\{s_n\}_{n=1}^\infty\sub \si \Bb$ satisfying 
$|s_n|\leq |f| $ for $n\in \N$ and $s_n\to f$ pointwise for $n \to \infty$. 
It follows from the Dominated Convergence Theorem for 
vector measures, \cite[Theorem 3.7(i)]{ORS}, that $\lim_{n\to \infty}s_n=f$ in the norm of $L^1(m_X)$. 
Consequently, 
we have both
\begin{equation}\label{3.18}
\lim_{n\to \infty} \int_A s_n\;d m_X=\int_Af \;d m_X \mbox{ for } A\in \Bb 
\mbox{ and } \lim_{n\to \infty}\Vert s_n\Vert_{m_X}=\Vert f \Vert_{m_X}
\end{equation}
because $\|\int_A(s_n-f)\,dm_X\|_X\le \|s_n-f\|_{m_X}$ for $A\in \Bb$, 
\cite[(3.21) on p.112]{ORS}, and
also $\big|\|s_n\|_{m_X}-\|f\|_{m_X}\big|\le \|s_n-f\|_{m_X}$ for $n\in\N$. 
Since the norms  $ \Vert \cdot \Vert_{[T,X]}$ and $\Vert \cdot \Vert_{m_X} $ 
coincide on $\si\Bb$ by (\ref{3.16}), the sequence $\{s_n\}_{n=1}^\infty$is also 
a Cauchy sequence in $[T,X]$ and hence, possesses a limit $h\in [T,X]$. 
So, $\lim_{n\to \infty} s_n=h$ also in the larger space $L^1$ 
(by definition $[T,X]\subseteq L^1$) because Step A of the proof of 
Lemma 4.4 in \cite{COR} shows that the natural inclusion $[T,X]\sub L^1$ is 
continuous. Thus, $\{s_n\}_{n=1}^\infty$ admits a subsequence converging a.e.\ 
to $h$. Therefore, $f=h$ a.e.\  on $(-1,1)$ and hence, $f\in [T,X]$ and 
$\lim_{n\to \infty} s_n=f$ in $[T,X]$. In particular, 
$\lim_{n\to \infty}\Vert s_n\Vert_{[T,X]}=\Vert f\Vert_{[T,X]}$. This, together with (\ref{3.16}) and (\ref{3.18}), yield
\begin{align*}
\Vert f\Vert_{[T,X]}=\lim_{n\to \infty}\Vert s_n\Vert_{[T,X]}=\lim_{n\to \infty} \Vert s_n\Vert_{m_X}=\Vert f\Vert_{m_X}.
\end{align*}
So, we have established the inclusion $L^1(m_X)\sub [T,X]$ as well as (\ref{3.14}).

Next, observe via (\ref{3.17}) with $A:= (-1,1)$, that the bounded linear operators 
$T\colon [T,X]\to X$ and $I_{m_X}\colon L^1(m_X) \to X$ coincide on the 
dense subspace $\si \Bb$ of $L^1(m_X)$. So, those two operators 
must coincide on $L^1(m_X)=\overline{\si\Bb}$. Hence, (\ref{3.15}) holds.

(ii) 
(a)$\Leftrightarrow$(b) This will follow once we establish the identity 
$L^1(m_X)=[T,X]_a$. Since $L^1(m_X)$ is a closed linear subspace of 
$[T,X]$ (by part (i)) and since $\si \Bb$ is dense in $L^1(m_X)$, we have 
$L^1(m_X)=\overline{\si\Bb}=[T,X]_b$. On the other hand, the continuous 
inclusion $X\sub [T,X]$ (see Lemma \ref{lemma2.6}(ii)) implies that 
$X_a\sub [T,X]_a$. In particular, $\si \Bb\sub [T,X]_a$ (cf. Lemma \ref{lemma2.3}), 
which yields $[T,X]_a=[T,X]_b$, \cite[Theorem I.3.13]{BS}. Thus we have $L^1(m_X)=[T,X]_a$.

(a)$\Rightarrow$(c) This follows from (\ref{3.11}) and (\ref{3.15}).

(c)$\Rightarrow$(d) Given $x^*\in X^*$, we shall show that the $\C$-valued, finitely 
additive set function $\langle \nu_f, x^*\rangle $ on $\Bb$ is $\sa$-additive. 
In view of Lemma \ref{lemma2.3} and \cite[Corollary I.4.2]{BS} the restriction 
$x^*|_{X_a}$ belongs to $X_a^*=X'$. So, there exists $g\in X'$ 
such that $g=x^*|_{X_a}$. Given $A\in \Bb$, 
as $\nu_f(A)=T(f\chi_A)\in X_a$ by assumption, Lemma \ref{lemma2.6}(iii) gives
\begin{align*}
\langle \nu_f(A), x^*\rangle &= \int_{-1}^1 T(f\chi_A) g \; d \mu= - \int_{-1}^1 f \chi_A T_{X'}(g)\; d \mu\\
&=-\int_A f T _{X'}(g)\; d \mu.
\end{align*}
So, $\langle \nu_f, x^*\rangle $ is $\sa$-additive because $f T_{X'}(g)\in L^1$; 
see (a)$\Leftrightarrow$(d) of Lemma \ref{lemma2.6}(i). Since $x^*\in X^*$ 
is arbitrary, (d) follows by the Orlicz-Pettis Theorem, \cite[Corollary I.4.4]{DU}.

(d)$\Rightarrow$(b) Let $A_n\downarrow \emptyset$ in $\Bb$. Then 
(d) implies that $\lim_{n\to \infty} \Vert\nu_f\Vert (A_n)= \Vert \nu_f\Vert(\emptyset)=0$, 
\cite[Proposition I.3]{Ri2}. So, it follows from Lemma \ref{lemma3.1} 
that $\Vert f \chi_{A_n}\Vert_{[T,X]}=\Vert \nu_f\Vert (A_n)\to 0$ as $n\to \infty$ and hence, $f\in [T,X]_a$.

(iii) We first claim that $f\in L^0$ belongs to $L^1_{w}(m_X)$ if and only if
\begin{equation}\label{3.19}
f T_{X'} (g) \in L^1, \hspace{5mm} g\in X'.
\end{equation}
Indeed, assume first that $f\in L^1_w(m_X)$. Then, given any $g\in X'\sub X^*$, 
the function $f$ is $\langle m_X, g\rangle$-integrable, which means that $f T_{X'}(g)\in L^1$ 
via Proposition \ref{proposition3.2}(ii). Conversely, assume that $f$ satisfies (\ref{3.19}). 
Given $x^*\in X^*$, find $g\in X'$ such that $x^*=g$ on $X_a$ as in the proof of 
(c)$\Rightarrow$(d) in part (ii). Since $\Rr(m_X)\sub X_a$ (see Proposition 
\ref{proposition3.2}(i)), we have, again by Proposition \ref{proposition3.2}(ii), that
\begin{align*}
\langle m_X, x^*\rangle(A) &=\langle m_X(A), x^*\rangle=\langle m_X(A), g\rangle=\langle m_X, g\rangle(A)\\
&= -\int_A T_{X'}(g) \; d \mu,\hspace{5mm} A\in \Bb.
\end{align*}
This and (\ref{3.19}) yield that $f$ is $\langle m_X, x^*\rangle$-integrable. 
Hence, $f\in L^1_w(m_X)$ as $x^*\in X^*$ is arbitrary. This completes the proof of the claim. 

We now obtain the inclusion $[T,X]\sub L^1_w(m_X)$ because every $f\in [T,X]$ 
satisfies (\ref{3.19}); see (a)$\Leftrightarrow$(d) in Lemma \ref{lemma2.6}(i). 
To prove the reverse inclusion, let $f\in L^1_w(m_X)$. For the Rybakov functional 
$g_0\in X'$ of $m_X$ as given in (\ref{3.4}), we may conclude that $f\in L^1$ 
because $|\langle m_X,g_0\rangle|=\mu$; see Proposition \ref{proposition3.2}(iii). 
This together with (\ref{3.19}) imply that $f\in [T,X]$; see (a)$\Leftrightarrow$(d) of Lemma \ref{lemma2.6}(i).

With the identity $L^1_w(m_X)=[T,X]$ established, let us verify that the corresponding 
norms are equal. Let $f\in L^1_w(m_X)$. Select a sequence $\{s_n\}_{n=1}^\infty \sub \si \Bb$ 
satisfying (\ref{2.5}). Since both $L^1_w(m_X)$ and $[T,X]$ have lattice norms and 
since $\Vert\, |s_n|\,\Vert_{m_X,w}=\Vert\, |s_n|\,\Vert_{m_X}=\Vert\, |s_n|\,\Vert_{[T,X]}$ 
for $n\in \N$ (see (\ref{3.16})), the Fatou property of both $L^1_w(m_X)$ and $[T,X]$ implies that
\begin{align*}
\Vert f  \Vert_{m_X,w}&=
\Vert\,| f |\,  \Vert_{m_X,w} = \sup_{n\in\N} \Vert \,|s_n|\,  \Vert_{m_X,w}= \sup_{n\in\N} \Vert \,|s_n|\,  \Vert_{[T,X]}\\
&=\Vert \,|f|\,  \Vert_{[T,X]}=\Vert f  \Vert_{[T,X]}.
\end{align*}
So, $L^1_w(m_X)$ and $[T,X]$ have equal norms. This completes the proof.
\end{proof}

Regarding the identity $[T,X]=L^1_w(m_X)$, the corresponding 
identity may fail to hold for kernel operators in general. See 
Example 3.4 of \cite{CR} which provides a counterexample 
for which the left-side is \emph{strictly smaller} than the right-side.

Let us present some immediate consequences of the previous theorem.

\begin{corollary}\label{corollary3.11}
Let $X$ be a r.i.\  space over $(-1,1)$ satisfying 
$0<\lbal_X\leq \ubal_X <1$. The following statements hold for 
the vector measure $m_X\colon \Bb\to X$.
\begin{enumerate}[(i)]
\item If $X $ has o.c. norm, then $X \sub L^1(m_X)=[T,X]$.
\item The B.f.s.\ $[T,X]$ has o.c. norm if and only if $L^1(m_X)=[T,X]$.
\item The natural inclusion $X\sub L^1(m_X)$ holds if and only if $T_X(X)\sub X_a$.
\item The natural inclusion $L^1(m_X)\sub X$ holds if and only if $L^1(m_X)=X_a$.
\item The equality $X=L^1(m_X)$ is valid if and only if 
$X$ has o.c.\ norm and equals $[T,X]$. In this case $L^1(m_X)=[T,X]$.
\item Suppose that $X=[T,X]$. Then $L^1(m_X)=X_a$. 
Hence, $X$ fails to have o.c.\ norm if and only if $L^1(m_X)\subsetneqq L^1_w(m_X).$
\item If $X$ satisfies $0<\lbal_X\leq \ubal_X <1/2$ or 
$1/2<\lbal_X\leq \ubal_X <1$, then the same conclusion as in (vi) holds.
\item Whenever $X=L^p$ for some $1<p<\infty$, we have 
\begin{equation*}
L^1(m_X) =L^1_w(m_X)=X=[T,X].
\end{equation*}
\item $L^1(m_X) =L^1(m_{X_a})$ with equal norms.
\end{enumerate}
\end{corollary}

\begin{proof}
(i) There are two methods to prove part (i).
\\
\emph{Method 1.} Since $X=X_a$, we have $X\sub L^1(m_X)$ by (\ref{3.13}). 
On the other hand, $L^1(m_X)=[T,X]$ is valid because $L^1(m_X)\sub[T,X]$ always
holds by (\ref{3.13}) whereas the reverse inclusion holds via Lemma 
\ref{lemma2.6}(i)  because of $X=X_a$ 
and the equivalence (a)$\Leftrightarrow$(c) in Theorem \ref{theorem3.10}(ii).
\\
\emph{Method 2.} The r.i.\ space $X$ does not contain an isomorphic copy of 
$c_0$ (see Lemma \ref{lemma2.1}), which (as noted prior to Remark \ref{remark3.8}) 
implies that $L^1(m_X)=L^1_w(m_X)$. Accordingly, $L^1(m_X)= L^1_w(m_X)=[T,X]\supseteq X$ 
by Theorem \ref{theorem3.10}(iii).

(ii) This follows from (a)$\Leftrightarrow$(b) in Theorem 
\ref{theorem3.10}(ii) because $[T,X]$ has o.c. norm if and only if $[T,X]_a=[T,X]$.

(iii) Assume that $X\sub L^1(m_X)$. Let $f\in X$. Then $T_X(f)= I_{m_X}(f)\in X_a$ 
by (\ref{3.11}) and (\ref{3.15}) and hence, the inclusion $T_X(X)\sub X_a$ holds. 
Conversely, suppose that $T_X(X)\sub X_a$. Then, given $f\in X\sub [T,X]$ and 
$A\in \Bb$, we have $T_X(f \chi_A)\in T_X(X)\sub X_a$. The equivalence 
(a)$\Leftrightarrow$(c) in Theorem \ref{theorem3.10}(ii) then ensures that 
$f\in L^1(m_X)$, which establishes the inclusion $X\sub L^1(m_X)$.

(iv) Assume that $L^1(m_X)\sub X$. Then the natural embedding $j\colon L^1(m_X)\to X$ is 
continuous . Indeed, consider any sequence $\{f_n\}_{n=1}^\infty\sub L^1(m_X)$ 
satisfying $\lim_{n\to \infty} f_n =0$ in $L^1(m_X)$ such that $\{j(f_n)\}_{n=1}^\infty$ 
converges to a function $h$ in $X$. Since each of the B.f.s.'  $L^1(m_X)$ and $X$ is 
continuously embedded into $[T,X]$ (via Theorem \ref{theorem3.10}(i) and 
Lemma \ref{lemma2.6}(ii), respectively) we have both $\lim_{n\to \infty}f_n=0 $ 
and $\lim_{n\to \infty}j(f_n) =h$ in $[T,X]$, which implies that $h=0$ a.e. So, $j(f)=0$. 
Hence, $j$ is continuous by the Closed Graph Theorem.

It is routine to verify the inclusion $(L^1(m_X))_a\sub X_a$ via the continuity of $j$. 
This implies that $L^1(m_X)=X_a$ because $L^1(m_X)=(L^1(m_X))_a$ and 
because $L^1(m_X)\supseteq X_a$ (see Lemma \ref{lemma3.6}). 
Conversely, if $L^1(m_X)=X_a$, then $L^1(m_X)\sub X$ as $X_a\sub X$.

(v) Assume that $X=L^1(m_X)$. In particular, $L^1(m_X)\sub X$ and so the natural 
inclusion is continuous (see the proof of part (iv)). By the Open Mapping 
Theorem and $L^1(m_X)=X$ it follows that $X$ and $L^1(m_X)$ are isomorphic B.f.s.'. 
So, $X$ necessarily has o.c. norm as $L^1(m_X)$ does. Moreover, by part (i) we have $X=L^1(m_X)= [T,X]$.

Conversely, assume that $X$ has o.c. norm and that $X=[T,X]$. Hence, $[T,X]$ also has o.c.\ norm. 
Then part (ii) gives the equality $X=L^1(m_X)$.

(vi) The assumption $X=[T,X]$ and continuity of the inclusion $X\sub [T,X]$ ensure that these two 
B.f.s.' are isomorphic (via the Open Mapping Theorem). Accordingly, $X_a=[T,X]_a$, which 
yields $X_a=L^1(m_X)$ because of $L^1(m_X)=[T,X]_a$ (see (\ref{3.13}) and 
equivalence (a)$\Leftrightarrow$(b) in Theorem \ref{theorem3.10}(ii)).

Next, suppose that $X$ fails to have o.c. norm. Then $L^1(m_X)\neq X$ 
by part (v). Moreover, $X_a=L^1(m_X)$ implies that $L^1(m_X)\subsetneqq X$. 
Since $[T,X]=L^1_w(m_X)$, by Theorem \ref{theorem3.10}(iii), and $X=[T,X]$ 
we can conclude that $L^1(m_X)\subsetneqq L^1_w(m_X)$.

Conversely, assume that $L^1(m_X)\subsetneqq L^1_w(m_X)$. 
Then Theorem \ref{theorem3.10}(iii) implies that 
$L^1(m_X)\subsetneqq L^1_w(m_X)= [T,X]=X$. 
Accordingly,  $L^1(m_X)\neq X$ and so $X$ fails to have o.c. norm by part (v).

(vii) Part (vi) is applicable because $X=[T,X]$; see Lemma \ref{lemma2.8}(i).

(viii) We know that $X=[T,X]$; see Lemma \ref{lemma2.8}(ii). 
Moreover, $X$ has o.c. norm. So, by part (v) we have $X=L^1(m_X)$. 
Accordingly, the equalities in (viii) hold as we also know that $L^1_w(m)=[T,X]$ 
by Theorem \ref{theorem3.10}(iii).

(ix) Recall that $\Rr(m_X)\subseteq X_a$; see Proposition \ref{proposition3.2}(i). 
So, $\int_As\,dm_{X_a}=\int_As\,dm_X$ for all $A\in \Bb$, $s\in\si\Bb$.
Moreover, the restriction to the closed subspace $X_a\subseteq X$ of 
each element of $X^*$ belongs to $X_a^*$ and conversely every element 
of $X_a^*$ has an extension to an element of $X^*$ with the same norm. 
From the definition of the spaces $L^1(m_X)$ and $L^1(m_{X_a})$ and their 
norms, together with the fact that $\si\Bb$ is dense in both spaces, it follows that 
$L^1(m_X)=L^1(m_{X_a})$ as vector spaces and with equal norms.
\end{proof}

\begin{remark}\label{remark3.12}
(a) Assume that the r.i.\  space $X$ has o.c. norm and satisfies $0<\lbal_X\leq \ubal_X <1$.

(a-1) The identity $L^1(m_X)=[T,X]$ in Corollary \ref{corollary3.11}(i) means, for 
the operator $T_X\colon X\to X$, that its optimal lattice domain $[T,X]$ (see 
Proposition \ref{proposition2.7}) and its o.c. optimal lattice domain $L^1(m_X)$ (see Lemma \ref{lemma3.6}) coincide.

(a-2) Parts (i) and (ii) of Corollary \ref{corollary3.11} yield that $[T,X]$ has o.c norm whenever $X$ has o.c. norm.

(a-3) Part (v) of Corollary \ref{corollary3.11} implies that $X=[T,X]$ if and only if $X=L^1(m_X)$. 
In other words, $T_X$ does not admit a continuous linear extension to any \emph{strictly larger} 
B.f.s.\ if and only if it does not admit a continuous linear extension to any \emph{strictly larger} B.f.s.\ with o.c. norm.

(b) Let $1<p<2$ or $2<p<\infty$. Then $X:= L^{p,r}$ for $1\leq r \leq\infty$  
satisfies $\lbal_X= \ubal_X =1/p$ and so $X=[T,X]$; see Lemma \ref{lemma2.8}(i).

(b-1) Let $1\leq r<\infty$. Then $X$ has o.c. norm and so part (vii) of 
Corollary \ref{corollary3.11} gives $X=X_a=L^1(m_X)$. Combined with part (v) 
of Corollary \ref{corollary3.11} we conclude that $X=L^1(m_X)=[T,X]$. Then Theorem \ref{theorem3.10}(iii) yields
\begin{align*}
X=L^1(m_X)=L^1_w(m_X)=[T,X].
\end{align*}

(b-2) Let $r:= \infty$. Then $X=L^{p,\infty}$ fails to have o.c. norm and
\begin{align*}
X_a=(L^{p,\infty})_a= (L^{p,\infty})_b= \left\{f\in L^{p,\infty}:\lim_{t\to 0+}t^{1/p}f^*(t)=0\right\};
\end{align*}
see \cite[Examples 3.4 and 3.10]{CR}. We then have
\begin{align*}
(L^{p,\infty})_a= L^1(m_X) \subsetneqq L^1_w(m_X)=[T,L^{p,\infty}]=L^{p,\infty}=X
\end{align*}
and hence, in particular, $L^1(m_X)\subsetneqq L^1_w(m_X)$. Indeed, 
$X=[T,X]$ and Corollary \ref{corollary3.11}(vi) yield $(L^{p,\infty})_a=L^1(m_X)$. 
Further, the proper inclusion $L^1(m_X)\subsetneqq L^1_w(m_X)$ is a 
consequence of Corollary \ref{corollary3.11}(vi) whereas $L^1_w(m_X)=[T,L^{p,\infty}]$ 
follows from Theorem \ref{theorem3.10}(iii).

(c) Let $p:= 2$. Then $\lbal_X= \ubal_X =1/2$ for $X:= L^{2,r} (1\leq r \leq \infty)$; see Section \ref{section2}.

(c-1) Consider the case when $1\leq r <2$ or $2<r< \infty$ 
(as the case $r:= 2$ gives $L^{2,2}=L^2$). Then $X= L^{2,r}$ has 
o.c.\  norm, so that $X=X_a$. Hence, Corollary \ref{corollary3.11}(i) gives
\begin{align*}
X\sub L^1(m_X)=[T,X].
\end{align*}

(c-2) Next consider the case when $r=\infty$. Then $X=L^{2,\infty}$ fails to have o.c. norm. Indeed,
\begin{align*}
(L^{2,\infty})_a=(L^{2,\infty})_b=\left\{f\in L^{2,\infty}: \lim_{t\to 0+}t^{1/2}f^*(t)=0\right\},
\end{align*}
again by \cite[Examples 3.4 and 3.10]{CR}. Then
\begin{align*}
(L^{2,\infty})_a\sub L^1(m_X)\sub [T,L^{2,\infty}],
\end{align*} 
by Theorem \ref{theorem3.10}(i). It is also the case that $L^{2,\infty} \sub [T,L^{2,\infty}]$; see Lemma \ref{lemma2.6}(ii).

\quad The function $1/w$ on $(-1,1)$ plays an important role in the study of the finite Hilbert transform. 
According to the discussion immediately after Lemma 2.1 in \cite{COR}, the space $L^{2,\infty}$ is the 
smallest r.i.\  space containing $1/w$. We note that $1/w\not\in(L^{2,\infty})_a$ because $(1/w)^*$ is 
the function $t\mapsto 2t^{-1/2}$ on $(0,2]$. It is interesting to know whether or not $1/w\in L^1(m_X)$? 
Since $1/w\in L^1_w(m_X)=[T,X]$, we see that $1/w\in L^1(m_X)$ if and only if 
$T(\chi_A (1/w))\in (L^{2,\infty})_a$ for all $A\in \Bb$; see (\ref{3.11}).

(d) For $X=L^p$ and $1<p<\infty$ most of the conclusion of part (viii) in 
Corollary \ref{corollary3.11} essentially occurs in \cite[Example 4.21]{ORS}. 
\end{remark}


\section{Integral representation of $T\colon [T,X]\to X$}\label{section4}
Integral representations for operators, when available, provide a useful 
method to investigate certain properties of the operator (e.g., compactness, 
complete continuity, etc.). Given a r.i.\  space $X$ over $(-1,1)$ satisfying 
$0<\lbal_X\leq \ubal_X <1$, recall from Theorem \ref{theorem3.10} that $L^1(m_X)$ 
is a closed subspace of $[T,X]$ and that $L^1(m_X)=[T,X]$ if and only if $[T,X]$ has 
o.c.\ norm. It can happen that $L^1(m_X)$ is strictly smaller than $[T,X]$ (see 
Remark \ref{remark3.12}(b-2)) in which case we \emph{cannot} represent 
$T\colon [T,X]\to X$ via integration with respect to the $X$\emph{-valued} 
vector measure $m_X$. Nevertheless, the purpose of this section is to 
present an integral representation of $T\colon [T,X]\to X$ of a ``different nature''. This requires an explanation.

When $L^1(m_X)$ is strictly smaller than $[T,X]=L^1_w(m_X)$, the idea is to enlarge the domain 
space $L^1(m_X)$ of $I_{m_X}$ to accommodate the extra functions from $[T,X]$ but the ``integrals'' 
of such functions should remain in $X$. One approach, as studied in \cite{CDR}, is to consider the 
bidual operator $I_{m_X}^{**}\colon L^1(m_X)^{**}\to X^{**}$, where $L^1(m_X)^{**}$ and $X^{**}$ are 
the biduals of the Banach spaces $L^1(m_X)$ and $X$, respectively. We may interpret $I_{m_X}^{**}$ as 
a ``sort of'' generalized integration operator. A relevant fact is that $[T,X]=L^1_w(m_X)$ is contained in 
$L^1(m_X)^{**}$, \cite[(10) on p.70]{CDR}. 
With this in mind, in order that the restriction operator $I_{m_X}^{**}|_{[T,X]}$ qualifies 
as a genuine generalized integral representation of $T\colon [T,X]\to X$, we would require the inclusion 
$I_{m_X}^{**}([T,X])\sub X$ to hold (as linear subspaces of $X^{**}$). But, it follows from 
\cite[Proposition 2.2]{CDR} that $I_{m_X}^{**}(L^1_w(m_X))\sub X$ if and only if $L^1(m_X)$ has the 
Fatou property. On the other hand, $L^1(m_X)$ has the Fatou property if and only if $L^1(m_X)=L^1_w(m_X)$; 
see Lemma \ref{lemma3.9}(ii). Accordingly, $I_{m_X}^{**}([T,X])\sub X$ holds if 
and only if $L^1(m_X)=L^1_w(m_X)=[T,X]$, in which case the operator $T\colon [T,X]\to X$ 
is exactly the 
integration operator $I_{m_X}$ (see Theorem \ref{theorem3.10}(i)). So, attempting to use 
$I_{m_X}^{**}$ as a possible $X$-valued extension of $I_{m_X}$ to a larger domain space 
which includes $[T,X]$ is not feasible. However, if we relax the requirement that the extension 
is $X$-valued, then the situation improves.

Accordingly, our strategy is to allow the possibility of enlarging the codomain space $X$ of the 
vector measure $m_X$. It turns out that this does enable us to obtain an integral representation 
of $T\colon [T,X]\to X$, whether $[T,X]$ has o.c. norm or not. We provide two natural examples of 
such extended codomain spaces (others are also possible). The first one is a space 
$\widehat{X}\supseteq X$ which is determined by $X$ and is minimal in the sense that, 
for \emph{every} $p\in (1,\infty)$ satisfying $X \sub L^p$, we necessarily have $\widehat{X}\sub L^p$. 
The second one is the space $L^0$. It is maximal and contains not only \emph{all} the r.i.\  
spaces over $(-1,1)$ but also $L^{1,\infty}$.

To define the Fr\'{e}chet  space $\widehat{X}\supseteq X$ we need some 
preparation. By a Fr\'{e}chet  space we mean a complete, metrizable, locally 
convex space. All Banach spaces are Fr\'{e}chet  spaces. For our r.i.\  space 
$X$ satisfying $0<\lbal_X\leq \ubal_X <1$, we know that $X\neq L^1$. 
Moreover, there always exists $1<p<\infty$ such that $X\sub L^p$ 
continuously by Lemma \ref{lemma2.2}. Note that the inclusion $X\sub L^p$ 
automatically guarantees that the inclusion is continuous via the Closed Graph Theorem. Define
\begin{align*}
J_X:= \{ p\in (1,\infty):X\sub L^p\},
\end{align*}
which is non-empty as noted above. Moreover, $J_X$ is a subinterval of $(1,\infty)$ 
because given any pair $p_1,p_2\in J_X$, the entire interval between $p_1$ and 
$p_2$ lies within $J_X$. Again by Lemma \ref{lemma2.2} we also have $L^q\sub X$ 
for some $1<q<\infty$. So, by defining 
\begin{align*}
p_X:= \sup J_X \leq q<\infty,
\end{align*}
we have 
\begin{align*}
X\sub L^p, \hspace{5mm} 1<p<p_X.
\end{align*}
Regarding $p_X$, there are two cases: either $p_X\in J_X$ or $p_X\not \in J_X$. 
Both cases occur. For example, $p_X=2$ if $X$ is any one of the spaces $L^{2,r}$ 
for $1\leq r\leq \infty$. However, if $1\leq r\leq 2$, then $p_X=2\in J_X=(1,2]$ 
whereas if $2<r\leq \infty$, then $p_X=2\not\in J_X=(1,2)$.

Let us return to the general r.i.\  space $X$. Define
\begin{align*}
\widehat{X}:= \bigcap_{p\in J_X} L^p.
\end{align*}
Equip $\widehat{X}$ with the locally convex Hausdorff topology induced by the 
norms $\Vert \cdot \Vert_{L^p}$ with $p$ varying throughout $J_X$. In practice, 
we do not require all
 $p\textrm{'}s$ in $J_X$ as explained below.

In the case when $p_X\in J_X$, it is clear that $\widehat{X}$ is normable and 
equals the reflexive Banach space $L^{p_X}$. In particular, $ (\widehat{X})^{*}=L^{(p_X)'}$

Next, assume that $p_X\not\in J_X$. In this case, $\widehat{X}$ is non-normable and is precisely 
the classical Fr\'{e}chet  space $L^{p_X-}$; see, for example, \cite{CDM}. To be more precise, take 
any increasing sequence $\{p_n\}_{n=1}^\infty$ in $J_X$ satisfying $\sup_{n\in \N} p_n=p_X$. 
Then the above defined locally convex Hausdorff topology on $\widehat{X}$ coincides with that 
generated by the sequence of norms $\Vert \cdot \Vert_{p_n}$ for $n\in \N$ and hence, is metrizable. 
As each $L^{p_n}$ for $n\in \N$ is a reflexive Banach space and 
$\widehat{X} =L^{p_X-}=\bigcap_{n=1}^\infty L^{p_n}$, it is known that $L^{p_X-}$ 
is a non-normable reflexive Fr\'{e}chet  space. Moreover, it is also known that 
\begin{equation}\label{4.1}
(\widehat{X})^{*}=\bigcup_{p\in J_X} L^{p'}=\bigcup_{n=1}^\infty L^{(p_n)'}.
\end{equation}
\quad In order to verify the continuous inclusion 
\begin{equation}\label{4.2}
[T,X]\sub \widehat{X},
\end{equation}
let us establish that 
\begin{equation}\label{4.3}
[T,X]\sub \bigcap_{p\in J_X} [T,L^p]=\bigcap_{p\in J_X}L^p=\widehat{X}.
\end{equation}
Given any $p\in J_X$, we have $X\sub L^p$ and so $[T,X]\sub [T,L^p]$; apply the equivalence 
(a)$\Leftrightarrow$(c) in Lemma \ref{lemma2.6}(i). Moreover, $[T,L^p]=L^p$ as isomorphic B.f.s.' by 
Corollary \ref{corollary3.11}(viii). So, (\ref{4.3}) holds. The continuity of the natural embedding from 
$[T,X]$ into $\widehat{X}$ follows from the fact that, for each $p\in J_X$, the inclusion $[T,X]\sub [T,L^p]$ 
is continuous by the definition of the respective norm in each space; see Lemma \ref{lemma2.6}(ii). 
This ensures that the inclusion (\ref{4.2}) is continuous because the topology on $\widehat{X}$ is 
determined by the family of norms $\{\Vert \cdot \Vert_{L^p}:p\in J_X\}$.

It is clear that $\widehat{X}$ is an order ideal of $L^0$ and that $\widehat{X}$ is a 
\emph{Fr\'{e}chet  function space} (briefly F.f.s) over $((-1,1),\Bb,\mu)$; see for example, 
\cite[\S 2.3]{B}, \cite{dCR} and the references therein. Observe that $X$ is continuously embedded 
into $\widehat{X}$ via (\ref{4.2}) together with the fact that the inclusion $X\sub [T,X]$ is continuous; 
see Lemma \ref{lemma2.6}(ii). This enables us to interpret $m_X\colon \Bb \to X$ as being 
$\widehat{X}$-valued. Moreover, the resulting set function $m_{\widehat{X}}\colon  \Bb\to \widehat{X}$ 
is $\sa$-additive, that is, $m_{\widehat{X}}$ is a F.f.s-valued vector measure. A set $A\in \Bb$ is said 
to be $m_{\widehat{X}}$-null if $m_{\widehat{X}}(B)=0$ for all $B\in \Bb$ contained in $A$.  
The claim is that the $m_{\widehat{X}}$-null sets and $\mu$-null sets coincide. Indeed, a set 
$A\in \Bb$ is $m_{\widehat{X}}$-null if and only if it is $m_{L^p}$-null for all $p\in J_X$ by the 
definition of $m_{\widehat{X}}$. The $m_{L^p}$-null and the $\mu$-null sets are the same for all 
$p\in J_X$ (by Proposition \ref{proposition3.2}(iv) with $X:= L^p$). The claim is thereby established.

Given a function $g\in (\widehat{X})^*$, let $\langle m_{\widehat{X}},g\rangle$ denote the $\C$-valued 
measure $A\mapsto \langle m_{\widehat{X}}(A),g\rangle$ on $\Bb$; see (\ref{4.1}). We say that 
$f\in L^0$ is scalarly $m_{\widehat{X}}$-integrable if $f$ is $\langle m_{\widehat{X}},g\rangle$-integrable 
for every $g\in (\widehat{X})^*$. Such a function $f$ is called $m_{\widehat{X}}$-integrable if, for every 
$A\in \Bb$, there exists a unique vector $\int_A f\; d m_{\widehat{X}}\in \widehat{X}$ satisfying 
$\langle \int_A f\; d m_{\widehat{X}},g\rangle = \int_Af \; d \langle m_{\widehat{X}},g\rangle$ for 
every $g\in (\widehat{X})^*$. Note that the reflexive F.f.s.\ $\widehat{X}$ cannot contain an 
isomorphic copy of $c_0$. This and \cite[Theorem 4]{Tu} enable us to apply \cite[Theorem 5.1]{L2} 
to conclude that every scalarly $m_{\widehat{X}}$-integrable is necessarily $m_{\widehat{X}}$-integrable. 
Let $L^1(m_{\widehat{X}})$ denote the vector space of all (equivalence classes of) 
$m_{\widehat{X}}$-integrable functions by identifying those functions which are $m_{\widehat{X}}$-a.e.\  
equal. It becomes a F.f.s.\ for the topology of uniform convergence of indefinite integrals; for the 
precise definition and details see, for example, \cite[\S 2.4]{B}, \cite{dCR}.

\begin{lemma}\label{lemma4.1}
Let $X$ be a r.i.\  space $X$ over $(-1,1)$ satisfying $0<\lbal_X\leq \ubal_X <1$.
\begin{enumerate}[(i)]
\item We have the continuous inclusions $X\sub [T,X]\sub \widehat{X}\sub L^1$.
\item The finite Hilbert transform $T$ maps $\widehat{X}$ into itself and the 
resulting linear operator $T\colon \widehat{X}\to \widehat{X}$ is continuous but, not compact.
\item The identity $\widehat{X}=L^1(m_{\widehat{X}})$ holds. In particular, 
$[T,X]\sub L^1(m_{\widehat{X}})$ with a continuous inclusion.
\item $[T,X]= \{f\in \widehat{X}: T(\chi_A f)\in X, \mbox{ for all } A\in \Bb\}$.
\end{enumerate}
\end{lemma}

\begin{proof}
(i) The continuity of the first two inclusions is already established; see (\ref{4.2}) 
and Lemma \ref{lemma2.6}(ii). The continuity of $\whX\sub L^1$ follows from the 
fact that $L^p\sub L^1$ continuously for all $p\geq 1$.

(ii) We have 
\begin{align*}
T(\whX)=T \bigg(\bigcap_{p\in J_X} L^p \bigg)\sub\bigcap_{p\in J_X} 
T\hspace{-0.7mm}\left(L^p \right)\sub\bigcap_{p\in J_X} L^p =\whX
\end{align*}
as we already know that $T(L^p)\sub L^p$ for all $p\in J_X$. Moreover, 
the operator $T\colon \whX\to\whX$ is continuous because, for each $p\in J_X$, we have
\begin{align*}
\Vert T(f)\Vert_{L^p}=\Vert T_{L^p}(f)\Vert_{L^p}\leq \Vert T_{L^p}\Vert_{\mbox{\tiny{op}}}
\Vert f\Vert_{L^p},\hspace{5mm} f\in \whX.
\end{align*}

Suppose that $T\colon \whX\to \whX$ is compact. According to the definition there exists a 
neighbourhood of $0$ in $\whX$ (which we may assume has the form 
$U:= \{f\in \whX:\Vert f\Vert_{L^p}\leq \eps\}$ for some $\eps>0$ and $p\in J_X$) 
such that $K:= \overline{T(U)}$ is compact in $\whX$. The continuity of the inclusion 
$\whX\sub L^p$ implies that $K$ is also compact in $L^p$, where we denote it by 
$K_p$. Observe that $U\sub V :=\{g\in L^p:\Vert g\Vert_{L^p}\leq \eps\}$ as 
$\whX\sub L^p$. Let $f\in V$ be given. Then there exists a sequence 
$\{s_n\}_{n=1}^\infty\sub \si \Bb$ such that $|s_n|\leq |f|$ for $n\in \N$ and 
$s_n\to f$ in $L^p$ for $n\to \infty$. By the continuity of $ T\colon L^p\to L^p$, 
also $T(s_n)\to T(f)$ in $L^p$. But, $\{s_n\}_{n=1}^\infty\sub U$ and so 
$\{T(s_n)\}_{n=1}^\infty\sub K\sub K_p$. By the compactness of $K_p$ it 
follows that $T(f)\in K_p$. Hence, $T(V)\sub K_p$. So, we have shown that 
$T\colon L^p\to L^p$ is compact, which is \emph{not} the case as its spectrum is 
uncountable, \cite[Theorem 13.9]{J}. Accordingly, $T\colon \whX\to \whX$ is not compact.

(iii) Corollary \ref{corollary3.11}(viii) gives the identities
\begin{equation}\label{4.4}
\whX= \bigcap_{p\in J_X} L^p=\bigcap_{p\in J_X} L^1(m_{L^p})=\bigcap_{p\in J_X} L^1_w(m_{L^p}).
\end{equation}
To prove the inclusion $\whX\sub L^1(m_{\whX})$, fix $f\in \whX$. Take any 
function $g\in (\whX)^*$. Via (\ref{4.1}) there exists $p\in J_X$ such that $g\in L^{p'}$. It is routine to check that
\begin{equation}\label{4.5}
\langle m_{\whX},g\rangle = \langle m_{L^p},g\rangle,
\end{equation} 
as $\C$-valued measures on  $\Bb$. Since $f\in L^p_w(m_{L^p})$, by (\ref{4.4}), 
we can conclude from (\ref{4.5}) that $f\in L^1(\langle m_{\whX},g\rangle )$. 
As $g\in (\whX)^*$ is arbitrary, the function $f$ is scalarly $m_{\whX}$-integrable 
and hence, also $m_{\whX}$-integrable; see the discussion prior to Lemma \ref{lemma4.1}. 
This establishes the inclusion $\whX\sub L^1(m_{\whX})$.

Conversely, let $f\in L^1(m_{\whX})$. Given any $p\in J_X$, it follows from (\ref{4.1}) 
and (\ref{4.5}) that $f\in L^1(\langle m_{L^p},g\rangle$ for every $g\in L^{p'}$. 
That is, $f\in L^1_w(m_{L^p})=L^p$ and hence, $f\in \whX$ because $p\in J_X$ 
is arbitrary. Thus, also $L^1(m_{\whX})\sub \whX$.
It can be checked that $ \whX=L^1(m_{ \whX})$ is also an isomorphism, as F.f.s.'

The continuity of $[T,X]\sub L^1(m_{\whX})$ now follows from $\whX=L^1(m_{\whX})$ and part (i).

(iv)  Clearly the right-side is contained in $[T,X]$ via the definition of $[T,X]$ and the fact that $\whX\sub L^1$.

To establish the reverse inclusion let $f\in [T,X]$. Then $f\in L^1$ and $T(\chi_A f)\in X$ 
for all $A\in \Bb$. Fix $p\in J_X$. Then $X\sub L^p$ and so $T(\chi_A f)\in L^p$ for all $A\in \Bb$. 
By definition of $[T,L^p]$ it follows that $f\in [T,L^p]=L^p$. Since $p\in J_X$ is arbitrary, 
we can conclude that $f\in \bigcap_{p\in J_X} L^p =:\whX$. 
\end{proof}

\begin{remark}\label{remark4.2}
(i) Let $X$ be a r.i.\  space $X$ over $(-1,1)$ satisfying $0<\lbal_X\leq \ubal_X <1$. 
Then the vector measure $m_{\whX}$ fails to have relatively compact range in $\whX$. 
If so, choose any $p\in J_X$. Continuity of the inclusion $\whX\sub L^p$ implies that 
$\{m_{\whX}(A):A\in \Bb\}=\{T(\chi_A):A\in \Bb\}$ is relatively compact in $L^p$. 
This contradicts Proposition \ref{proposition3.2}(vi) (with $X:=L^p$ there). 
Hence, $\Rr(m_{\whX})$ is \emph{not} relatively compact in $\whX$.

(ii) If $p\not \in J_X$, then the inclusion $[T,X]\sub \whX=L^1(m_{\whX})$ (cf. (\ref{4.2})) is \emph{proper}. 
Indeed, in the event of equality, the continuous inclusion of the Banach space $[T,X]$ into the Fr\'{e}chet  
space $\whX$ would be surjective. By the Open Mapping Theorem it would be an isomorphism 
which contradicts the fact that $\whX$ is non-normable. The analogous argument (still for $p_X\not\in J_X$) 
shows that the inclusion $\whX\sub L^1$ is also proper. The same is true if $p_X\in J_X$ as then 
$\whX=L^{p_X}$ with $p_X>1$ and hence, $\whX$ cannot equal $L^1$.

(iii) There is an alternative way to see that the inclusion $\whX\sub L^1(m_{\whX})$ is continuous by 
applying an optimal domain theorem for continuous linear operators on F.f.s.'. Indeed, it is easy to verify 
that $\whX$ has a $\sa$-Lebesgue topology (i.e., $\lim_{n\to \infty} f_n$ in $\whX$ whenever 
$f_n\downarrow 0$ in the order of $\whX$), \cite[Example 2.3.1]{B}. The vector measure 
$m_{\whX}\colon \Bb\to \whX$ satisfies $m_{\whX}(A)=T(\chi_A)$ for $A\in \Bb$, where we 
may consider $T$ as the continuous linear operator $T\colon \whX\to \whX$; see 
Lemma \ref{lemma4.1}(ii). As observed earlier, the $m_{\whX}$-null and the $\mu$-null sets 
coincide; in other words, $T$ is $\mu$\emph{-determined}. These facts allow us to apply the 
Optimal Domain Theorem, \cite[Theorem 3.3.1]{B}, to deduce that $\whX\sub L^1(m_{\whX})$. 
Actually, Theorem 3.3.1 in \cite{B} asserts more. Namely, that $L^1(m_{\whX})$ is the \emph{largest} 
F.f.s.\ with a $\sa$-Lebesgue topology into which $\whX$ is continuously embedded and to 
which $T\colon \whX\to\whX$ admits an $\whX$-valued continuous linear extension, 
realized as the integration operator $I_{m_{\whX}}\colon L^1(m_{\whX})\to \whX$, that is, 
\begin{align*}
I_{m_{\whX}}(f)= \int_{-1}^1 f\; d m_{\whX}=T(f), \hspace{5mm} f\in \whX.
\end{align*}
\end{remark}

We are ready to present an integral representation for $T\colon [T,X]\to X$ via the 
vector measure $m_{\whX}\colon \Bb\to \whX$.

\begin{proposition}\label{proposition4.3}
Let $X$ be a r.i.\  space $X$ over $(-1,1)$ satisfying $0<\lbal_X\leq \ubal_X <1$. 
Then every function $f\in [T,X]\sub L^1(m_{\whX})$ satisfies 
\begin{align*}
\int_Af \; d m_{\whX}= T(f \chi_A)\in X, \hspace{5mm} A\in \Bb.
\end{align*}
Consequently, the linear operator $T\colon [T,X]\to X$ admits the $X$-valued integral representation given by 
\begin{equation}\label{4.6}
T(f)= \int_{-1}^1 f\; d m_{\whX}\in X,\hspace{5mm} f\in [T,X].
\end{equation}
\end{proposition}

\begin{proof}
Let $f\in [T,X]$. Then $f\in L^1(m_{\whX})=\whX$ by Lemma \ref{lemma4.1}(iii). 
Let $g\in (\whX)^*=\bigcup_{p\in J_X} L^{p'}$ (see (\ref{4.1})). Select any 
$p\in J_X$ such that $g\in L^{p'}$. Given $A\in \Bb$, we apply both (\ref{3.15}) of 
Theorem \ref{theorem3.10}(i), with $X:= L^p$, and (\ref{4.5}) to obtain that 
\begin{align*}
\langle\int_A f \; d m_{\whX},g \rangle &=\int_A f \; d\langle m_{\whX},g \rangle
=\int_A f \; d\langle m_{L^p},g \rangle=\langle\int_A f \; d m_{L^p},g \rangle\\
&=\langle T_{L^p}(f\chi_A),g \rangle=\langle T(f\chi_A),g \rangle,
\end{align*}
where we used the fact that $f\in [T,X]$ implies (by definition of $[T,X]$) that 
$T(f\chi_A)\in X\sub L^p$. As $g\in (\whX)^*$ is arbitrary and $X\sub \whX$, 
we can conclude that $\int_A f \; d m_{\whX}=T(f\chi_A)$ as an equality in $\whX$. 
But, actually $T(f\chi_A)\in X\sub \whX$ and so $\int_A f\;d m_{\whX}=T(f\chi_A)\in X$ 
for all $A\in \Bb$ and $f\in [T,X]$. Apply this to the case when $A:= (-1,1)$ to 
deduce the integral representation (\ref{4.6}).
\end{proof}

Next, consider $L^0$ as an extended codomain space for $m_X$. It is maximal in 
the sense that it contains \emph{all} the r.i.\  spaces over $(-1,1)$ as well as 
$L^{1,\infty}$ which is the codomain space of $T\colon L^1\to L^{1,\infty}$. 
Recall that $L^0$ is is a topological vector space for the topology of convergence 
in measure, which is defined by the $F$-norm
\begin{align*}
f \mapsto \int_{-1}^1 \frac{|f|}{1+|f|}\; d\mu, \hspace{5mm} f\in L^0.
\end{align*}
As $L^0$ is complete, it is an $F$-space by definition. For a theory of integration 
with respect to an $F$-space-valued measure see, for example, \cite{CD}, and the references therein.

Let $X$ be a r.i.\  space over $(-1,1)$ satisfying $0<\lbal_X\leq \ubal_X <1$. 
The natural embedding from $X$ into $L^0$ is continuous; see \cite[Theorem I.1.4]{BS}. 
So, the vector measure $m_X$ considered as $L^0$-valued is 
$\sa$-additive. We denote this vector measure by $m_X^0\colon \Bb\to L^0$. 
The integral $\int_A s \; d m_X^0$ of each $s\in \si \Bb$ over a 
set $A\in \Bb$ is defined in the canonical way.  In view 
of \cite[Proposition 3.1]{CD}, a $\Bb$-measurable function 
$f\colon (-1,1)\to \C$ is said to be $m_X^0$-integrable if there 
exists a sequence $\{s_n\}_{n=1}^\infty$ in $\si \Bb$ satisfying 
$\lim_{n\to \infty} s_n=f$ pointwise and such that, for every 
$A\in \Bb$, the sequence $\{\int_A s_n\; d m_X^0\}_{n=1}^\infty$ 
converges in $L^0$. In this case the integral $\int_A f \, d m_X^0$ 
of $f$ over $A$ is defined as the limit $\lim_{n\to \infty} \int_A s_n \;d m_X^0$ in 
$L^0$ (which is independent of the choice of such a sequence $\{s_n\}_{n=1}^\infty$).

We conclude this section with a second integral 
representation of the operator $T\colon [T,X]\to X$, namely via $m_X^0$.

\begin{proposition}\label{proposition4.4}
Let $X$ be a r.i.\  space over $(-1,1)$ satisfying $0<\lbal_X\leq \ubal_X <1$.  
\begin{enumerate}
\item[(i)] Every $f\in L^1$ is $m_X^0$-integrable and $\int_Af\,dm_X^0=T(f\chi_A)\in L^0$ for
$A\in\Bb$.
\item[(ii)] If $f\in[T,X]\subseteq L^1$, then  $\int_Af \; d m_X^0\in X$ for every 
$A\in \Bb$. So, the linear operator $T\colon [T,X]\to X$ admits the integral representation 
\begin{equation}\label{4.7}
T(f)= \int_{-1}^1 f \; d m_X^0\in X, \hspace{5mm} f\in [T,X].
\end{equation}
\end{enumerate}
\end{proposition}

\begin{proof}
(i) By linearity it suffices to consider $0\le f\in L^1$. Select a sequence $\{s_n\}_{n=1}^\infty
\subseteq\si\Bb$ such that $0\le s_n\uparrow f$ pointwise on $(-1,1)$. 
By the Monotone Convergence Theorem $\|f-s_n\|_{L^1}\to0$ for $n\to\infty$. 
For $A\in\Bb$ it follows that $s_n\chi_A\to f\chi_A$ in $L^1$ for $n\to\infty$. 
Since $T\colon L^1\to L^0$ is continuous, we can conclude that $T(s_n\chi_A)\to T(f\chi_A)$ in $L^0$.
But, $T(s_n\chi_A)=\int_A s_n\,dm_X^0$ for $n\in\N$ and so 
$\int_A s_n\,dm_X^0\to T(f\chi_A)$ in $L^0$. 
By definition it follows that $f$ is $m_X^0$-integrable and $\int_A f\,m_X^0= T(f\chi_A)$ for each $A\in\Bb$.

(ii) Suppose  $f\in[T,X]$. Then $T(f\chi_A)\in X$ for every $A\in\Bb$ and so, from
part (i), $\int_Af\,dm_X^0\in X$. The choice $A=(-1,1)$ yields \eqref{4.7}.
\end{proof}

Since $T\colon L^1\to L^{1,\infty}$ is continuous, as is the natural inclusion $L^{1,\infty}\subseteq L^0$,
the proof of part (i) of Proposition \ref{proposition4.4} shows that actually 
$\int_Af\,dm_X^0\in L^{1,\infty}$ for every $f\in L^1$ and $A\in\Bb$. As noted in Section 2, the containment $[T,X]\subseteq L^1$ is proper.
Moreover, $[T,X]=L^1_w(m_X)$ (see Theorem 3.1(iii)) 
shows that every scalarly $m_X$-integrable function
is also $m_X^0$-integrable.



\begin{thebibliography}{99}

\bibitem{BS}
C. Bennett, R. Sharpley, \emph{Interpolation of Operators}, Academic
Press,  Boston, 1988.

\bibitem{B}
B. Blaimer \emph{Optimal Domain and Integral Extension of 
Operators Acting in Fr\'echet Function Spaces}, Logos Verlag, 
Berlin; Katholische Univ. Eichst\"att-Ingolstadt (Ph D Thesis), 2017.
Also available at https://zenodo.org/record/1087454 .

\bibitem{CD}
G. P. Curbera, O. Delgado, \emph{Optimal domains for $L^0$-valued 
operators via stochastic measures}, Positivity \textbf{11} (2007), 399--416.


\bibitem{CDM}
J. M. F. Castillo, J. C. Diaz, J. Motos, \emph{On the Fr\'{e}chet  
space $L_{p-}$}, Manuscripta Math. \textbf{96} (1998), 219--230.


\bibitem{CDR}
G. P. Curbera, O. Delgado, W. J. Ricker, \emph{Vector measures: 
Where are their integrals?}, Positivity \textbf{13} (2009), 61--87.

\bibitem{COR}
G. P. Curbera, S. Okada, W. J. Ricker, \emph{ Inversion and
extension of the finite Hilbert transform on $(-1,1)$},
Ann. Mat. Pura Appl. \textbf{198} (2019), 1835--1860.

\bibitem{CR-N}
G. P. Curbera, W. J. Ricker, \emph{Optimal domains for kernel
operators via interpolation}, Math. Nachr., \textbf{244} (2002),
47--63.

\bibitem{CR} G. P. Curbera, W. J. Ricker, \emph{Banach lattices with the Fatou property
and optimal domain of kernel operators}, Indag. Math., N.S.,
 \textbf{17} (2006),  187--204.

\bibitem{CR-Survey}G.P. Curbera, W.J. Ricker,  \emph{Vector measures, integration and
applications}. In: Positivity, Eds. K. Boulabiar, G. Buskes and  A.
Triki, Trends in Mathematics, Birkh\"auser Verlag,
Basel-Berlin-Boston, 2007, pp. 127--160.

\bibitem{dCR}
R. del Campo, W. J. Ricker , \emph{The Fatou completion of a Fr\'{e}chet  
function space and applications}, J. Aust. Math. Soc., \textbf{88} (2010), 49--60.

\bibitem{DU}
J. Diestel, J.J. Uhl, Jr., \emph{Vector Measures}, Math. Surveys \textbf{15},
Amer. Math. Soc., Providence, R.I., 1977.

\bibitem{J} K. J\"orgens
\textit{Linear Integral Operators}, (English transl.) Pitman, Boston, 1982.

\bibitem{L1}
D.R. Lewis,  \emph{Integration with respect to vector measures},
Pacific J. Math.,  \textbf{33} (1970), 157--165.

\bibitem{L2} D.R. Lewis,  \emph{On integrability and summability in vector spaces}, Illinois J.
Math., \textbf{16} (1972), 294--307.

\bibitem{LZII}
J. Lindenstrauss,  L. Tzafriri, \emph{Classical Banach Spaces Vol.
II}, Springer-Ver\-lag, Berlin, 1979.

\bibitem{LZ} W.A.J. Luxemburg,  A.C.  Zaanen, \textit{Riesz Spaces
I}, North-Holland, Amsterdam, 1971.

\bibitem{OE}
S. Okada, D. Elliot, \emph{The finite Hilbert transform in
$\mathcal{L}^2$}, Math. Nachr., \textbf{153} (1991),  43--56.

\bibitem{ORS}
S. Okada, W.J. Ricker, E.A. S\'anchez-P\'erez, \emph{Optimal Domain
and Integral Extension of Operators: Acting in Function Spaces},
Operator Theory Advances and Applications \textbf{180},
Birkh\"auser,  Berlin, 2008.

\bibitem{St}
G.F. Stefansson,  $L\sb 1$ of a vector measure, \textit{Le
Matematiche {\rm(}Catania{\rm)}} \textbf{48} (1993),  219--234.

\bibitem{Ri}
W.J. Ricker, \textit{Separability of the $L^1$-space of a vector measure}, 
Glasgow Math. J., \textbf{34} (1992), 1-9.

\bibitem{Ri2}
W.J. Ricker, \textit{Operator Algebras Generated by Commuting 
Projections: A Vector Measure Approach}, LNM 1711, Springer, Berlin Heidelberg, 1999.

\bibitem{T}
F.G. Tricomi, \emph{Integral Equations}, Interscience, New York,
1957.

\bibitem{Tu}
Yu. B. Tumarkin, \emph{On locally convex spaces with basis}, 
Dokl. Akad. Nauk SSSR \textbf{195} (1970), 1278--1281 (in Russian): 
English transl.: Soviet Math. Dokl. \textbf{11} (1970), 1672--1675.


\bibitem{vD} D. van Dulst, \textit{Characterizations of Banach Spaces 
not Containing $\ell^1$}, CWI Tract No.59, Centrum voor Wiskunde en Informatica, Amsterdam, 1989.

\bibitem{Z} A.C. Zaanen,  \textit{Integration}, 2nd rev. ed., North-Holland, Amsterdam, 1967.

\bibitem{ZII}
A.C. Zaanen, \emph{Riesz Spaces II}, North-Holland, Amsterdam, 1983.

\end{thebibliography}
\end{document}